\documentclass{article}

\usepackage{color}
\usepackage{amsmath}
\usepackage{amsfonts}
\usepackage{amssymb}
\usepackage{amsbsy}
\usepackage{bm}
\usepackage{graphicx}
\usepackage{hyperref}
\usepackage{nicefrac}

\newcommand{\R}{\mathbb{R}}
\newcommand{\uv}{\ensuremath{\hat{u}}} 
 
\newcommand{\vv}{\ensuremath{\hat{v}}} 
\newcommand{\mv}{\ensuremath{\hat{m}}} 
\newcommand{\hf}{\nicefrac{1}{2}}

\usepackage{authblk}
\title{Convexity splitting in a phase field model for surface diffusion}

\author[1]{Rainer~Backofen}
\author[2]{Steven~M.~Wise}
\author[1]{Marco~Salvalaglio}
\author[1]{Axel~Voigt}
\affil[1]{
  Department of Mathematics,
  Technische Universit\"at Dresden
  01062 Dresden, Germany}

\affil[2]{
  Department of Mathematics,
 The University of Tennessee, 
 Knoxville, TN 37996, USA 
}

\begin{document}







\maketitle

\begin{abstract}
{Convexity splitting like schemes with improved accuracy are proposed for a phase field model for surface diffusion. The schemes are developed to enable large scale simulations in three spatial dimensions describing experimentally observed solid state dewetting phenomena. We introduce a first and a second order unconditionally energy stable scheme and carefully elaborate the loss in accuracy associated with large time steps in such schemes. We then present a family of Rosenbrock convex splitting schemes. We show the existence of a maximal numerical timestep and demonstrate the increase of this maximal numerical time step by at least one order of magnitude using a Rosenbrock method. This scheme is used to study the effect of contact angle on solid state dewetting phenomena.}
\end{abstract}

\section{Introduction}
If an energy ${\mathcal{E}}$ can be written as the difference of two convex energies ${\mathcal{E}}_c$ and ${\mathcal{E}}_e$, ${\mathcal{E}} =  {\mathcal{E}}_c - {\mathcal{E}}_e$, then the time discretization
\begin{equation*}
\frac{u^{n+1} - u^n}{\tau_n} = - \nabla_{\mathcal H} {\mathcal{E}}_c [u^{n+1}] - \nabla_{\mathcal H} {\mathcal{E}}_e [u^n]
\end{equation*}
of the gradient flow
\begin{equation*}
u_t = - \nabla_{\mathcal H} {\mathcal{E}}[u]
\end{equation*}
is energy stable. That is, it satisfies ${\mathcal{E}}[u^{n+1}] \leq {\mathcal{E}}[u^n]$ for all time steps $n$. Here, $\tau_n = t^{n+1} - t^n$ is the discrete time step width, $u^n$ denotes the time-discrete approximation of $u(t^n)$, $\nabla_{\mathcal H} {\mathcal{E}}$ denotes the gradient of an energy ${\mathcal{E}}$ with respect to the inner product of a Hilbert space $\mathcal{H}$, defined by $(\nabla_{\mathcal H} {\mathcal{E}}[u], \theta)_{\mathcal H}$ $= \delta_u {\mathcal{E}} [u](\theta)$ for all $\theta \in {\mathcal H}$, where the right-hand side is the Gateaux derivative of ${\mathcal{E}}$ in a test function $\theta$. Typical choices are ${\mathcal H} = L^2$, for non-conserved flows,  and ${\mathcal H} = H^{-1}$, for conserved flows. This convexity splitting idea is often attributed to Eyre \cite{Eyre_MRS_1998} and has become popular as a simple and efficient discretization scheme for various evolution problems with a gradient flow structure, see e.g. \cite{Chengetal_JCP_2008,Wiseetal_SIAMJNA_2009,Gomezetal_JCP_2011,Elseyetal_ESAIM_2013,Vignal2015,Shin2016,Shin2017}. Some of these schemes are shown to be unconditionally energy stable, unconditionally solvable and converge optimally in the energy norm. However, it has also been shown that convexity splitting approaches can lead to large errors \cite{Chengetal_JCP_2008,Christliebetal_CMS_2013,Elseyetal_ESAIM_2013,Glasneretal_JCP_2016}. We will elaborate on this issue and propose a convexity splitting approach for a phase field model for surface diffusion with improved accuracy.

The model to be considered reads \cite{Raetzetal_JCP_2006}
\begin{eqnarray}
\label{eq1}
\partial_t u = \nabla \cdot \mathbf{j}, \qquad
\mathbf{j} = \frac{1}{\epsilon} M(u) \nabla \mu, \qquad
g(u)\mu = \frac{1}{\epsilon} B^\prime(u) - \epsilon \Delta u,
\end{eqnarray}
in $\Omega \times(0,\infty)$ with $\Omega \subset \R^d$ with $d = 2,3$. The variable $u$ denotes an order parameter for the phases of the system, such that, for example, $u = 1$ represents a solid phase, and $u = 0$ represents a liquid phase.  The variable $\mu$ is the chemical potential. We consider the initial condition $u(\mathbf{x},0) = u_0(\mathbf{x})$ and boundary conditions, i.e. $\mathbf{n} \cdot \nabla \mu = \mathbf{n} \cdot \nabla u = 0$, where $\mathbf{n}$ is the outward unit normal to $\partial \Omega$. $B(u) = 18 u^2(1- u)^2$ is a double-well potential, $M(u) = 2 B(u)$ a mobility function, $g(u) = 30 u^2(1-u)^2$ an enhancing function, and $\epsilon$ relates to the thickness of the transition region between the two phases $u = 1$ and $u = 0$. The model formally converges for $\epsilon \to 0$ to Mullins sharp interface model for surface diffusion \cite{Mullins_JAP_1957}, see, in particular \cite{Raetzetal_JCP_2006,Guggenbergeretal_PRE_2008,Voigt_APL_2016}. For the aformentioned result to hold the fourth order polynomial in $u$ in the mobility function $M(u)$ is essential \cite{Voigt_APL_2016}. As recently shown \cite{Daietal_MMS_2014,Leeetal_APL_2015,Leeetal_SIAMJAM_2016} occasionally used second order polynomials in $u$ in the mobility function $M(u)$, see e.g. \cite{Bhateetal_JAP_2000,Wiseetal_JCP_2006,Torabietal_PRSA_2009,Jiangetal_AM_2012}, do actually not converge to surface diffusion if $\epsilon \to 0$. Heuristic arguments and matched asymptotic analysis lead to the presence of an additional bulk diffusion term, which might alter the long time behavior. This is not the case for the originally proposed phase field approximation \cite{Cahnetal_EJAM_1996}, which uses a double-obstacle potential instead of the double-well potential $B(u)$. The enhancing function $g(u)$ does not alter the matched asymptotic analysis. Such a function is commonly used in classical phase field models for solidification \cite{Karmaetal_PRE_1998} to ensures better, through not necessarily higher convergence in $\epsilon$. Herein we demonstrate its advantage to increasing the accuracy, especially for larger values of $\epsilon$. 

For $M(u) = g(u) = 1$ the classical Cahn-Hilliard equation \cite{Cahnetal_JCP_1958} is recovered, which formally converges for $\epsilon \to 0$ to the Mullins-Sekerka problem \cite{Mullinsetal_JAP_1963}, see \cite{Pego_PRS_1989}. The Cahn-Hilliard equation can be written in the abstract framework as the gradient flow $u_t = - \nabla_{H^{-1}} {\mathcal{E}}[u]$ as the $H^{-1}$, with respect to the energy
\begin{eqnarray}
{\mathcal{E}}(u) = \int_\Omega \left( \frac{\epsilon}{2} |\nabla u|^2 + \frac{1}{\epsilon} B(u) \right) \; d \mathbf{x}.
\end{eqnarray}
Convexity splitting for the Cahn-Hilliard equation has been considered in e.g. \cite{Eyre_MRS_1998,Bertozzietal_IEEE_2007,Glasneretal_JCP_2016}. In this work we start with the canonical nonlinear convex splitting $B(u)=B_c(u)-B_e(u)$, where
\begin{eqnarray}
\label{eq4}
B_c(u)=B(u)+ \alpha  \left(u-\frac{1}{2}\right)^2, \qquad
B_e(u)=\alpha \left(u-\frac{1}{2}\right)^2.
\end{eqnarray}
For $\alpha \geq 9$ the two function $B_c(u)$ and $B_e(u)$ are convex and thus also the energies ${\mathcal{E}}_c = \int_\Omega \left(\frac{\epsilon}{2}|\nabla u|^2 + \frac{1}{\epsilon} B_c(u)\right) d \mathbf{x}$ and ${\mathcal{E}}_e = \int_\Omega \left(\frac{1}{\epsilon} B_e(u)\right) d \mathbf{x}$. The resulting convex splitting scheme is unconditionally energy stable, unconditionally solvable and converges optimally in the energy norm \cite{Glasneretal_JCP_2016}. We will adapt this scheme and use it for eq. \eqref{eq1} with the above considered functional forms for $M(u)$ and $g(u)$. For $g(u) = 1$ the analytical results for the Cahn-Hilliard equation can be adapted to show energy stability properties also for the degenerate model \cite{elliott96a}. In this paper we recall a first and a second order unconditionally energy stable scheme for this case. But, if $g(u)$ in non-constant, we are unable to rigorously demonstrate the desired properties of the schemes, the model does not fall into the considered class of a gradient flow. It does not even have an energetic-variational structure. However, our goal is to construct practical and stable numerical schemes that enable large scale simulations in three spatial dimensions.  Our methods are already used in \cite{Naffoutietal_SA_2017} to provide predictive, validated simulation results for complex solid-state dewetting scenarios of ultra-thin silicon films.\footnote{For an introduction to solid-state dewetting see the review \cite{Thompson_ARMR_2012}.} To enable these simulations the improved accuracy with $g(u)$ is absolutely essential. We therefore adapt the proposed scheme, consider a Rosenbrock time discretization to increase the accuracy and apply it to the general case.

The paper is organized as follows: In Section \ref{sec2} we describe the numerical approach in detail. We propose a first and a second order scheme, for which we proof unconditional energy stability for the case $g(u) = 1$. A modified linear first order scheme and a new Rosenbrock time stepping scheme is then introduced for the more general case. In Section \ref{sec3} we analyze the schemes with respect to accuracy, solvability and efficiency. We consider the example of a retracting step in two space dimensions to find optimal parameters, which are used in Section \ref{sec4} for large scale simulations for solid-state dewetting in three space dimensions. We further discuss an outlook to more realistic modeling approaches including  the incorporation of vapor-substrate and film-substrate interfacial energies and illustrate the framework required to treat anisotropic energies. In Section \ref{sec5} we draw conclusions.
 
 \section{Discretization schemes}
\label{sec2}

Eq. \eqref{eq1} with the convex splitting in eq. \eqref{eq4}  can be written as a system of two second order partial differential equations
\begin{eqnarray}
\partial_t u = \nabla \cdot \left(\frac{1}{\epsilon} M(u) \nabla \mu\right), \label{eqn-SD-1} \qquad
 g(u) \mu = - \epsilon \Delta u + \frac{1}{\epsilon} B^\prime_c(u) - \frac{1}{\epsilon} B^\prime_e(u). 
\end{eqnarray} 
We observe that this equation does not have an energy dissipation structure, unless $g(u)$ is a constant function. On the other hand, the sharp interface law to which it converges is a gradient flow, and it is thus reasonable to expect stable dynamics. The question is: \emph{With what energy do we measure the stability of solutions to \eqref{eqn-SD-1}?} We currently cannot answer this question in a systematic, quantitative way and, therefore we recall the analysis only for the case $g(u) = 1$. Strictly speaking even with $g(u) = 1$ eq. \eqref{eqn-SD-1} is not a gradient flow, because of the presence of the mobility function $M(u)$. On the other hand, it does have a clear energy-variational structure, and this is all that is needed to discuss the issue of energy stability. With $M(u)$ the rate of dissipation is $d_t	\mathcal{E}[u] = -\int_\Omega \frac{1}{\epsilon}M(u)|\nabla \mu|^2d\mathbf{x}$ and thus ${\mathcal{E}}[u^{n+1}] \leq {\mathcal{E}}[u^n]$.

\subsection{A first order unconditionally energy stable scheme}
	
Here we review what might be considered a standard framework for constructing first order convex splitting schemes. Consider the convex splitting scheme
	\begin{align}
\frac{u^{n+1}-u^n}{\tau_n} = & \ \nabla \cdot \left( \frac{1}{\epsilon} M(u^n) \nabla \mu^{n+1}\right),
	\label{scheme-1-alt}
	\\
\mu^{n+1} = & \ - \epsilon \Delta u^{n+1} + \frac{1}{\epsilon} B^\prime_c(u^{n+1}) - \frac{1}{\epsilon} B^\prime_e(u^n).
	\label{scheme-2-alt}
	\end{align}
Since $B_c$ and $B_e$ are convex, it follows that $(u^{n+1}-u^n)B'_c(u^{n+1}) = B_c(u^{n+1}) -B_c(u^n) + R_c(u^n,u^{n+1})$,
$-(u^{n+1}-u^n)B'_e(u^n) =  -B_e(u^{n+1}) +B_e(u^n) + R_e(u^n,u^{n+1})$, where $R_c, R_e\ge 0$.  Testing \eqref{scheme-1-alt} with $\mu^{n+1}$ and \eqref{scheme-2-alt} with $u^{n+1}-u^n$, and using the boundary conditions, we obtain 
	\begin{align*}
\left(u^{n+1}-u^n, \mu^{n+1}\right)  = & -\frac{\tau_n}{\epsilon}\left(M (u^n)\nabla \mu^{n+1}, \nabla \mu^{n+1} \right)
	\\
\left(\mu^{n+1}, u^{n+1}-u^n \right) = & \ \epsilon\left(\nabla u^{n+1},\nabla\left(u^{n+1}-u^n\right)\right) + \frac{1}{\epsilon} \left( B^\prime_c(u^{n+1}),u^{n+1}-u^n \right) 
	\\
&  - \frac{1}{\epsilon}\left( B^\prime_e(u^n),u^{n+1}-u^n\right) 
 	\\
= &  \ \mathcal{E}[u^{n+1}] - \mathcal{E}[u^n]  + \int_\Omega  \frac{1}{\epsilon} \left(  R_c(u^n,u^{n+1}) +  R_e(u^n,u^{n+1}) \right) d\mathbf{x}
	\\
& + \frac{\epsilon}{2}\left( \nabla(u^{n+1}-u^n),\nabla(u^{n+1}-u^n)\right).
	\end{align*}
Combining both, we obtain the following energy stability result
	\begin{equation}
\mathcal{E}[u^{n+1}] + \frac{\tau_n}{\epsilon}\left(M (u^n)\nabla \mu^{n+1}, \nabla \mu^{n+1} \right) + \mathcal{R}^n_\epsilon = 	\mathcal{E}[u^n] 
	\label{energy-stability-1}
	\end{equation}
with 
\[
\mathcal{R}^n_\epsilon =  \int_\Omega  \frac{1}{\epsilon} \left(  R_c(u^n,u^{n+1}) +  R_e(u^n,u^{n+1}) \right) d\mathbf{x} +   \frac{\epsilon}{2}\left(\nabla(u^{n+1}-u^n),\nabla(u^{n+1}-u^n)\right)
\]
 the non-negative energy dissipation term.

\subsection{A second order unconditionally energy stable scheme}
	
Following \cite{guo16,yan18}, we can use the energetic-variational structure of eq. \eqref{eqn-SD-1}  and the convex splitting methodology to devise second order accurate (in time) energy stable schemes. Consider, in particular
	\begin{align}
\frac{u^{n+1}-u^n}{\tau} = & \ \nabla\cdot\left( \frac{1}{\epsilon}M (u^{n+\hf})\nabla \mu^{n+\hf}\right)
	\label{scheme-2nd-1-alt}
	\\
\mu^{n+\hf} = &  - \frac{3\epsilon}{4} \Delta u^{n+1} - \frac{\epsilon}{4} \Delta u^{n-1} 
 + \frac{1}{\epsilon} \frac{B_c(u^{n+1})- B_c(u^n)}{u^{n+1}-u^n}  - \frac{1}{\epsilon} B_e'({u}^{n+\hf}) 
	\label{scheme-2nd-2-alt}
	\end{align}
where ${u}^{n+\hf}$ is obtained via extrapolation
${u}^{n+\hf} = \frac{3}{2}u^n - \frac{1}{2} u^{n-1}$.
Considerations similar to those for the first order scheme are taken in the construction of this second order scheme. The convex terms are treated implicitly, and the concave term and the mobility are handled via extrapolation. One can prove that the scheme is unconditionally energy stable in the sense that 
	\begin{align}
\mathcal{F}[u^{n+1},u^n] \, + &  \, \frac{\tau}{\epsilon}\left(M({u}^{n+\hf})\nabla \mu^{n+\hf}, \nabla \mu^{n+\hf} \right) + \mathcal{R}^{n+\hf}_\epsilon = \  \mathcal{F}[u^n,u^{n-1}],
	\label{energy-stability-2}
	\end{align}
where $\mathcal{R}^{n+\hf}_\epsilon$ is a non-negative remainder term and $\mathcal{F}$ is a numerical energy defined as
	\begin{equation}
\mathcal{F}[u ,v] :=  \mathcal{E}[u] + \frac{1}{16\epsilon}\left\|(u -v)\right\|_{L^2}^2 +\frac{\varepsilon}{8}\| \nabla (u - v) \|_{L^2}^2.
	\label{numerical-energy}
	\end{equation}
This type of energy stability is sometimes called weak energy stability, because it involves a modification of the energy $\mathcal{E}[u]$. See \cite{diegel2016} for more details, including an optimal order convergence analysis for the case that $M(u) = 1$. One drawback of this scheme is that it is not a one-step scheme, which makes it particularly difficult to adapt the time step size $\tau>0$. It will therefore not be used for the full problem, as time-adaptivity will be essential for the large scale simulation, see Section \ref{sec4}.

\subsection{Convex splitting like schemes for the full problem}

Since our primary goal is to conduct efficient and accurate simulations over large time and space scales, employing temporal and spatial adaptivity will be critical. In order to construct practical, stable, and high order time stepping strategies in this setting, we turn to semi-implicit Runge-Kutta schemes. Such schemes often have excellent stability properties, can be of arbitrarily high order, and since they are one-step methods, can handle time adaptivity easily. Semi-implicit Runge-Kutta schemes that respect the convex-concave structure of the energy have been recently introduced in the literature~\cite{Shin2017}. These schemes, in their current incarnation, can be shown to be rigorously unconditionally energy stable and solvable, but only when $M(u)$ and $g(u)$ are a constant functions. Our schemes, based on the Rosenbrock framework, will be shown to be stable and accurate when $M(u)$ and $g(u)$ are non-constant. Moreover, they possess a natural error indicator that is useful for time step adaptivity, as we will show.
We adapt the proposed first order scheme and approximate the non-linear term to obtain a modified linear system. 

Rewrite eq. \eqref{eqn-SD-1} for $\uv = (u,\mu)$ as 
\begin{eqnarray}
H \partial_t \uv = F(\uv), \quad \mbox{with} \quad 
F(\uv) = F_c(\uv) +  F_e(\uv)
\end{eqnarray}
and
\begin{eqnarray*}
H=\begin{bmatrix} 0 & 0 \\ 1 & 0 \end{bmatrix}, \quad
F_c(\uv)=
\begin{bmatrix}
g(u) \mu + \epsilon \Delta u - \frac{1}{\epsilon} B'_c(u)\\
\nabla \cdot \left( M(u) \nabla \mu \right)
\end{bmatrix} \quad\text{and} \quad
F_e(\uv)=
\begin{bmatrix}
\frac{1}{\epsilon} B'_e(u)\\
0
\end{bmatrix} 
\end{eqnarray*}
we will consider a Taylor expansion to treat $F_c$ semi-implicitly. However, $M(u)$ and $g(u)$ are treated explicitly. We thus obtain the semi-implicit convex splitting like scheme  
\begin{eqnarray}
\label{eq:sics}
\frac{1}{\tau} H  \uv^{n + 1}  -F^*_{c,\uv}(\uv^{n}) \uv^{n + 1} = \frac{1}{\tau} H \uv^{n} + F_{c}(\uv^n)-F^*_{c,\uv}(\uv^{n}) \uv^{n}+ F_{e}(\uv^n)
\end{eqnarray}
with
\begin{eqnarray}
\label{eq:sie-jac}
F^*_{c,\uv}(\uv^{n})\uv^{n+1}=
\begin{bmatrix}
g(u^{n}) \mu^{n+1}+\epsilon \Delta u^{n+1} -\frac{1}{\epsilon} B''_c(u^{n})u^{n+1}\\
\nabla \cdot \left( M(u^{n}) \nabla \mu^{n+1} \right) 
\end{bmatrix}.
\end{eqnarray}
It can be considered as variant of the convex splitting  scheme~\eqref{scheme-1-alt} -- \eqref{scheme-2-alt}, with a Newton linearization to treat the cubic term and the inclusion of the function $g(u)$. Since $g(u)$, $B_c''(u)$, and $M(u)$ are positive functions, the solvability of this scheme can be established, at least in the spatially discrete case and under certain reasonable assumptions. The stability of the scheme is, however, unknown, as we have already stated.

We now use this first order scheme to build higher order Rosenbrock-Wanner schemes, which are semi-implicit Runge-Kutta schemes, which do not require iterative Newton steps, see \cite{Lan99} for a review. With these schemes one can achieve higher order methods for stiff problems by working the Jacobian matrix of $F$, in our case only $F^*_{c,\uv}$, or approximations of it, into the integration formula. For a general introduction in the context of ordinary differential equations we refer to the textbooks \cite{HairerWanner,DeuflhardBornemann}. Rosenbrock-Wanner schemes are defined by the recursive update 
\begin{equation}
  \uv^{n + 1} = \uv^{n} + \sum_{i = 1}^{s} m_{i} \uv^{n}_i\,
\end{equation}
where $\uv_i = (u_i,\mu_i)$ are the implicitly defined Runge-Kutta stages, $m_{i} \in \mathbb{R}$ are coefficients determining the order and accuracy of the method, and $s$ the stage order of the scheme. We now have to solve a system of $s$ partial differential equations for the unknown stages $\uv^{n}_1, \ldots, \uv^{n}_s$ in each time step
\begin{equation}
  \label{eq:rb-stage}
  \frac{1}{\tau} \frac{1}{\gamma} H \uv^{n}_i - F_{c,\uv}^*(\uv^{n}) \uv^{n}_i = 
   \frac{1}{\tau} H \sum_{j=1}^{i - 1} c_{ij} \uv^{n}_j + F_c(\vv^{n}_i) + F_e(\vv^{n}_i)
\end{equation}
with 
\begin{eqnarray}
\vv^{n}_i = \uv^{n} + \sum_{j = 1}^{i - 1}a_{ij}\uv^{n}_{j}
\end{eqnarray} 
and $\gamma$, $c_{ij}$ and $a_{ij}$ coefficients defined by the particular Rosenbrock method, see Tables \ref{tab1} and \ref{tab2}. We will consider a two-stage method (ROS2) and a three-stage method (ROS34PW2), see \cite{Verwereta_sjsc_1999,Rang_BIT_2005,Lan99}, and for both methods an approximation of the Jacobian
\begin{eqnarray}
\label{eq:ros-jac}
F^*_{c,\uv}(\uv^{n}) \uv^{n}_i = \begin{bmatrix}
g(u^n) \mu^{n}_i + g'(u^n) \mu^{n}  u^{n}_i + \epsilon \Delta u^{n}_i  - \frac{1}{\epsilon}B''_c(u_{n})u^{n}_i \\
\nabla \cdot \left(M'(u^{n}) u^{n}_i \nabla \mu^{n} \right) + \nabla \cdot \left( M(u^{n}) \nabla \mu^{n}_i \right) 
\end{bmatrix} 
\end{eqnarray}
where now the nonlinear terms $M(u)$ and $g(u)$ are treated semi-implicitly. However, numerical tests have shown that also the use of the simpler approximation eq. \eqref{eq:sie-jac} leads to similar results.

Our proposed Rosenbrock-Wanner schemes and the schemes from \cite{Shin2017} may be viewed in the same general semi-implicit Runge-Kutta framework, and the stability analysis may be related.

\begin{table}[ht]
\begin{center}
\begin{tabular}{|l|l|}
\hline 
$\gamma = 1.707106781186547$ & $c_{11} = \gamma$ \\ 
\cline{1-1}
$a_{21} = 0.5857864376269050$ & $c_{21} = -1.171572875253810$ \\
$a_{22} = 1.0$ & $c_{22} = -\gamma$ \\
 \hline
$m_1 = 0.8786796564403575$ & $\hat{m}_1 = 0.5857864376269050$ \\ 
$m_2 = 0.2928932188134525$ & $\hat{m}_2 = 0.0$ \\ 
\hline 
\end{tabular} 
\end{center}
\caption{A set of coefficients for the ROS2 Rosenbrock scheme. All coefficients not given explicitly are set to zero.}\label{tab1}
\end{table}

\begin{table}[ht]
\begin{center}
\begin{tabular}{|l|l|}
\hline 
$\gamma = 0.43586652150845$ & $c_{11} =  \gamma$ \\ 
\cline{1-1}
$a_{21} = 2$ & $c_{21} = -4.58856072055809$ \\
$a_{22} = 0.871733043016918$ & $c_{22} = -\gamma$ \\
$a_{31} = 1.41921731745576$ & $c_{31} = -4.18476048231916$ \\
$a_{32} = -0.25923221167297$ & $c_{32} = 0.285192017355496$ \\
$a_{33} = 0.731579957788852$ & $c_{33} = -0.413333376233886$\\
$a_{41} = 4.18476048231916$ & $c_{41} = -6.36817920012836$\\
$a_{42} = -0.285192017355496$ & $c_{42} = -6.79562094446684$\\
$a_{43} = 2.29428036027904$ & $c_{43} = 2.87009860433106$\\
$a_{44} = 1.0$ & $c_{44} = 0$\\
 \hline
$m_1 = 4.18476048231916$ & $\hat{m}_1 = 3.90701053467119$ \\ 
$m_2 = -0.285192017355496$ & $\hat{m}_2 = 1.1180478778205$ \\ 
$m_3 = 2.29428036027904$ & $\hat{m}_3 = 0.521650232611491$ \\ 
$m_4 = 1.0$ & $\hat{m}_4 = 0.5$ \\ 
\hline 
\end{tabular} 
\end{center}
\caption{A set of coefficients for the ROS34PW2 Rosenbrock scheme. All coefficients not given explicitly are set to zero.}\label{tab2}
\end{table}

The standard semi-implicit and Rosenbrock schemes, without convexity splitting, are obtained by setting $\alpha=0$ in eq. \eqref{eq4}. However, the resulting schemes are not solvable with standard iterative solvers for reasonable timesteps, which is the reason to introduce the convexity splitting schemes.

To discretize in space we use globally continuous, piecewise linear Lagrange finite elements and a conforming triangulation of the domain $\Omega$. To assemble and solve the resulting systems we use the FEM-toolbox AMDiS \cite{Veyetal_CVS_2007,Witkowskietal_ACM_2015}. The mesh is adaptively refined within the diffuse interface to ensure a minimal resolution with respect to the interface width $\epsilon$. As linear solver we used a BiCGStab(l) method, with $l=2$, and a  block Jacobi preconditioner with ILU factorization.

\section{Results}
\label{sec3}

\begin{figure}[htb]
\noindent
\centering
\begin{tabular}{c}
\includegraphics*[angle = -0, width = 0.85 \textwidth]{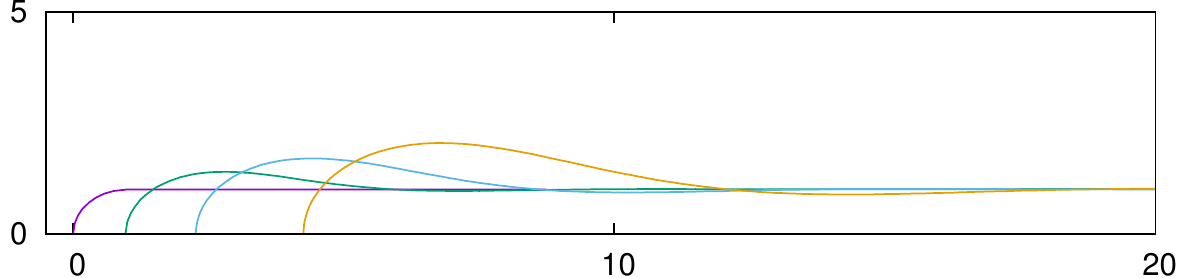} 
\end{tabular}
\begin{center}
\begin{minipage}{0.9\textwidth}
\caption[short figure description]{
Reference solution for a retracting step at $t = 0$, $2$, $9$ and $30$ from left to right. 
\label{partfig1}
}
\end{minipage}
\end{center}
\end{figure}

We consider an example in two space dimensions, which mimics solid-state dewetting. Fig. \ref{partfig1} shows the simulation results for a retracting step. The initial setup is a step with a small aspect ratio (height/length) on a substrate with $90^\circ$ contact angle. This configuration is out of equilibrium and will evolve towards a minimal energy state. A hill is formed followed by a small valley. As the step retracts, th ehill growth and the valley deepens, which eventually will lead to a splitting into well separated parts, each converging to a hemisphere. For larger aspect ratios this splitting can be prevented and the shape evolves towards its equilibrium shape which is given by the Winterbottom construction \cite{Winterbottom_AM_1967}, in our case again a hemisphere, as we consider a contact angle of $90^\circ$, see Section \ref{sec4} for more details. However, in order to analyze the numerical schemes we are here only interested in the early stage of evolution, the retracting step. We will use the tip position to validate our numerical approach. 

\subsection{Justification of modeling}
\begin{figure}[htb]
\noindent
\centering
\begin{tabular}{c}
\includegraphics*[angle = -0, width = 0.45 \textwidth ]{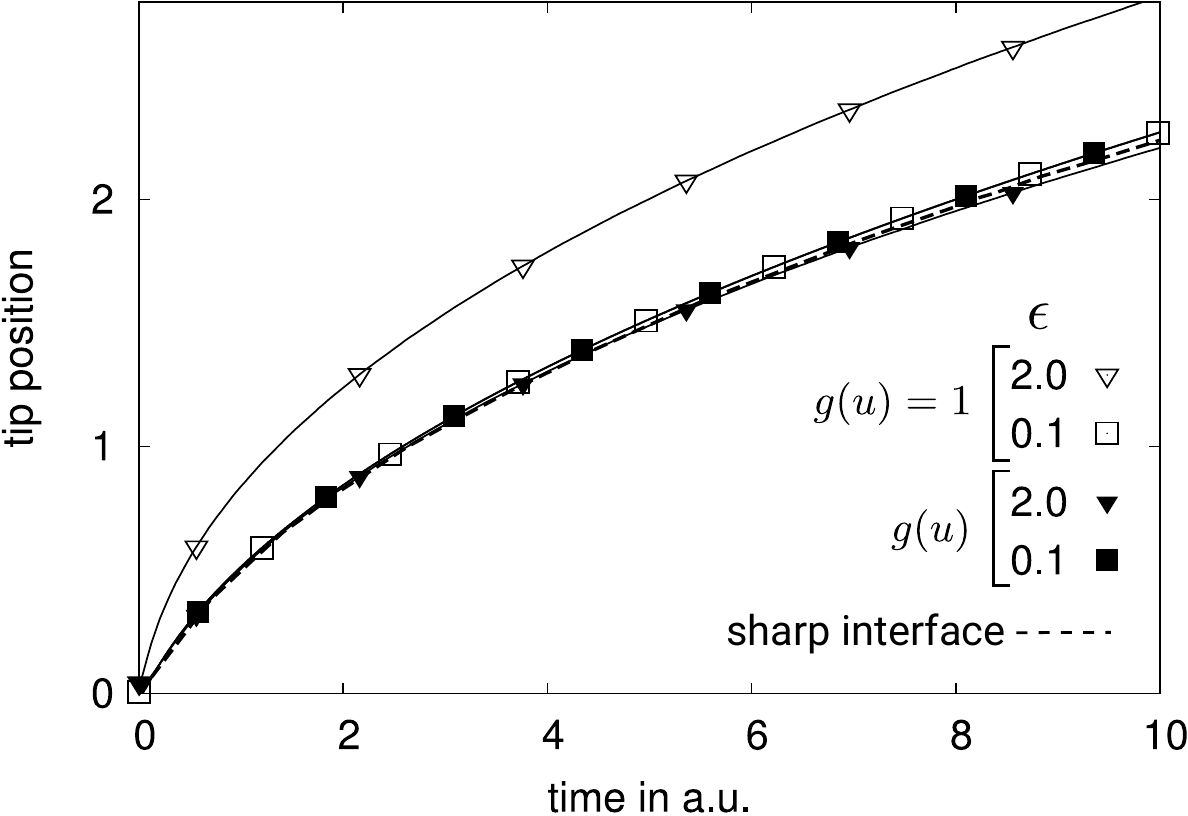}
\end{tabular}
\begin{center}
\begin{minipage}{0.9\textwidth}
\caption[short figure description]{
Tip position for simulations with and without $g(u)$ for two different values of $\epsilon$ in comparison with a reference solution. 
\label{figmodel}
}
\end{minipage}
\end{center}
\end{figure}

We first demonstrate the modeling advantage by comparing solutions with and without $g(u)$ for various $\epsilon$ with the reference solution, which here is the corresponding sharp interface solution for surface diffusion.  For numerical treatment of the sharp interface problem we refer e.g. to \cite{Baenschetal_JCP_2005,Hausseretal_IFB_2005,Salvalaglioetal_PRB_2016,Baoetal_JCP_2017}. Fig. \ref{figmodel} shows the comparison for various $\epsilon$, for the semi-implicit scheme with and without $g(u)$, $\tau = \epsilon 10^{-3} $, $\alpha = 0$ and a direct solver. A scheme, which leads to accurate results but is impractical for simulations in three space dimensions. The reference solution shows a $t^{1/2}$ scaling for the tip position. While this can only be reproduced for small values of $\epsilon$ if $g(u) = 1$, the behavior is also found for large values of $\epsilon$ if $g(u)$ is considered. 

\subsection{Comparision of convex-splitting schemes}
\begin{figure}[htb]
\noindent
\centering
\begin{tabular}{cc}
\includegraphics*[angle = -0, width = 0.45 \textwidth ]{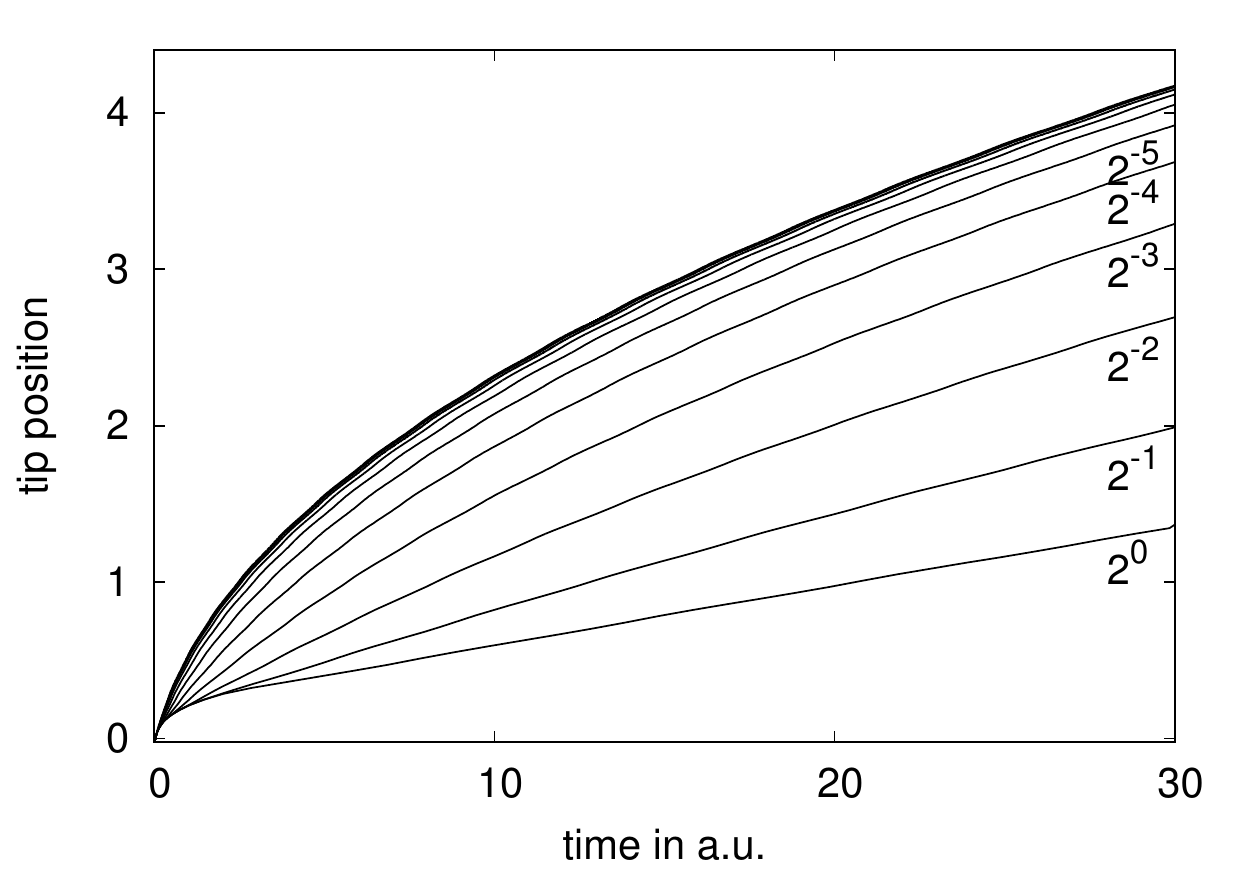} &
\includegraphics*[angle = -0, width = 0.45 \textwidth]{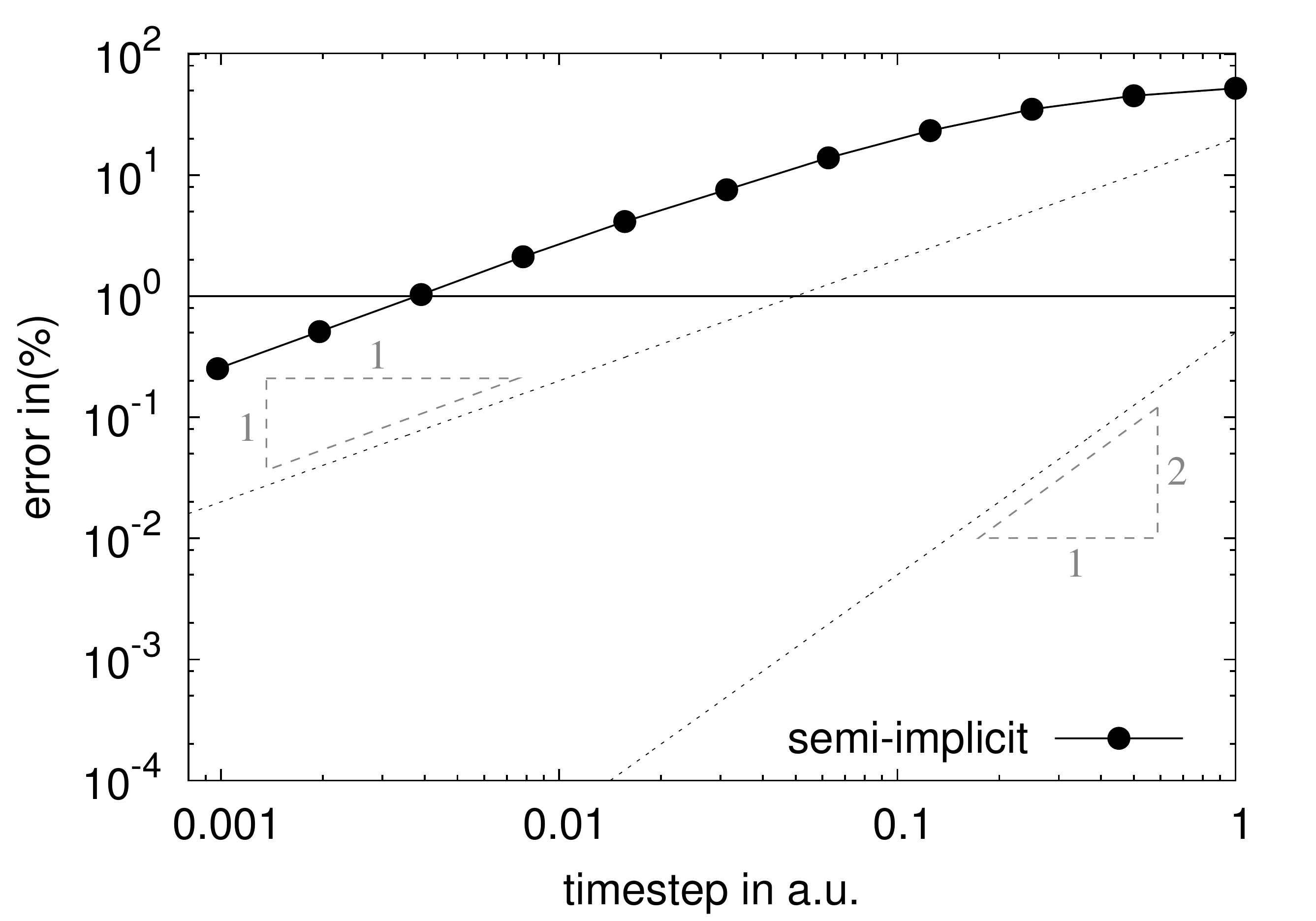} \\
\end{tabular}
\begin{center}
\begin{minipage}{0.9\textwidth}
\caption[short figure description]{
Semi implicit convexity splitting scheme. (left) Tip position over time for various timesteps together with the reference solution. The numbers indicate the used timesteps. (right) Deviation from the reference solution over timesteps.
\label{partfig2}
}
\end{minipage}
\end{center}
\end{figure}

We now compare the three proposed convex-splitting-like schemes. For the reference solutions we now consider the corresponding scheme with a small timestep $\tau = 10^{-4}$, $\alpha = 0$ and a direct solver.  All simulations start with a small time step $\tau = 10^{-4}$. After the initial phase at $t = 0.05$ the time step is gradually increased until it reaches the final time step $\tau$, which is reported in the following. We consider $\alpha = 9$, for which the resulting linear systems for each scheme and each $\tau$ can be solved by the mentioned iterative solver. Fig. \ref{partfig2} shows the tip position over time for various timesteps and the corresponding error for the semi implicit convexity splitting scheme. The correct qualitative behavior is only achieved by the semi implicit convexity splitting scheme for $\tau < 3\cdot10^{-2}$. To achieve a quantitative error, which is below 1\% even requires $\tau < 4\cdot10^{-3}$. Such large errors have also been reported for other convex-splitting schemes \cite{Chengetal_JCP_2008,Elseyetal_ESAIM_2013}. The results in Fig. \ref{partfig2} further indicate first order convergence in $\tau$. 

\begin{figure}[htb]
\noindent \centering
\begin{tabular}{cc}
\includegraphics*[angle = -0, width = 0.45 \textwidth ]{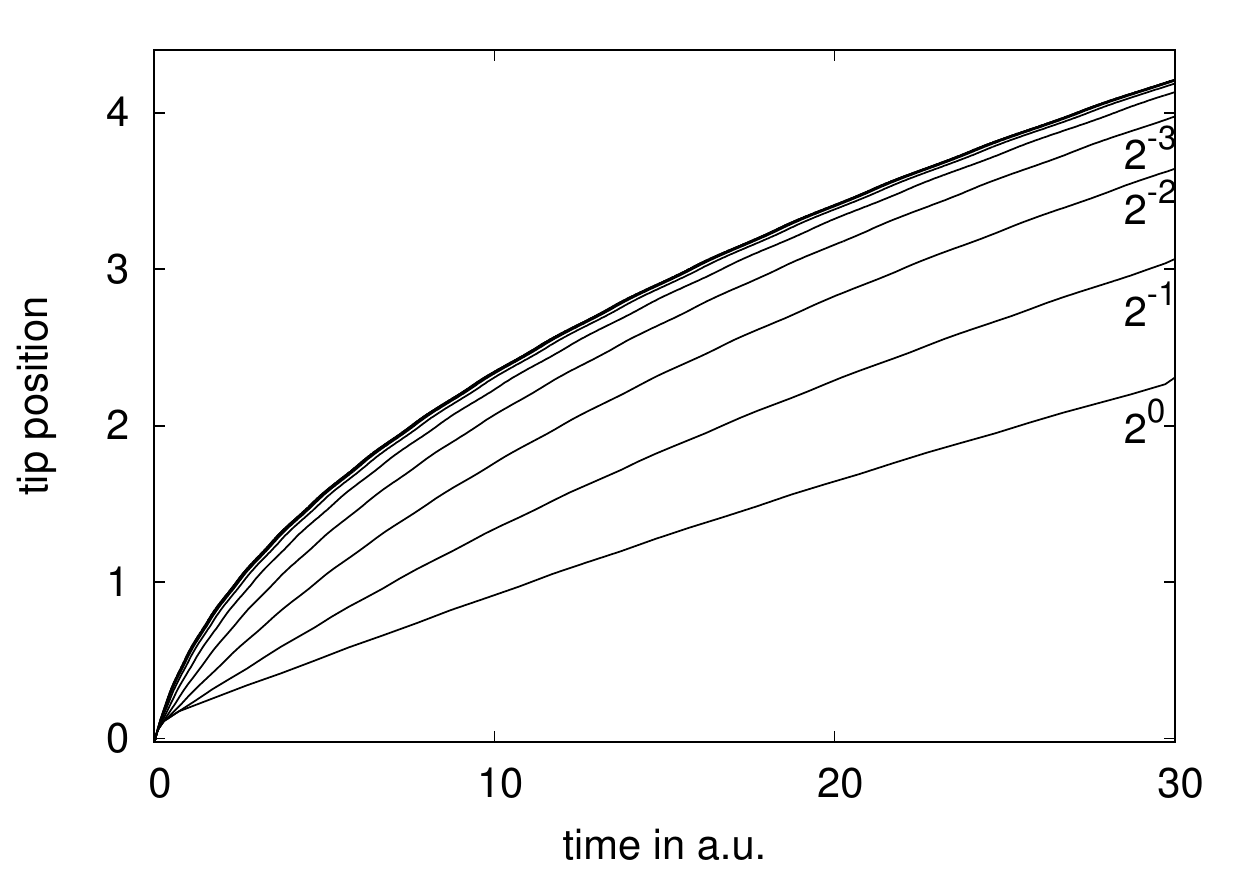} &
\includegraphics*[angle = -0, width = 0.45 \textwidth]{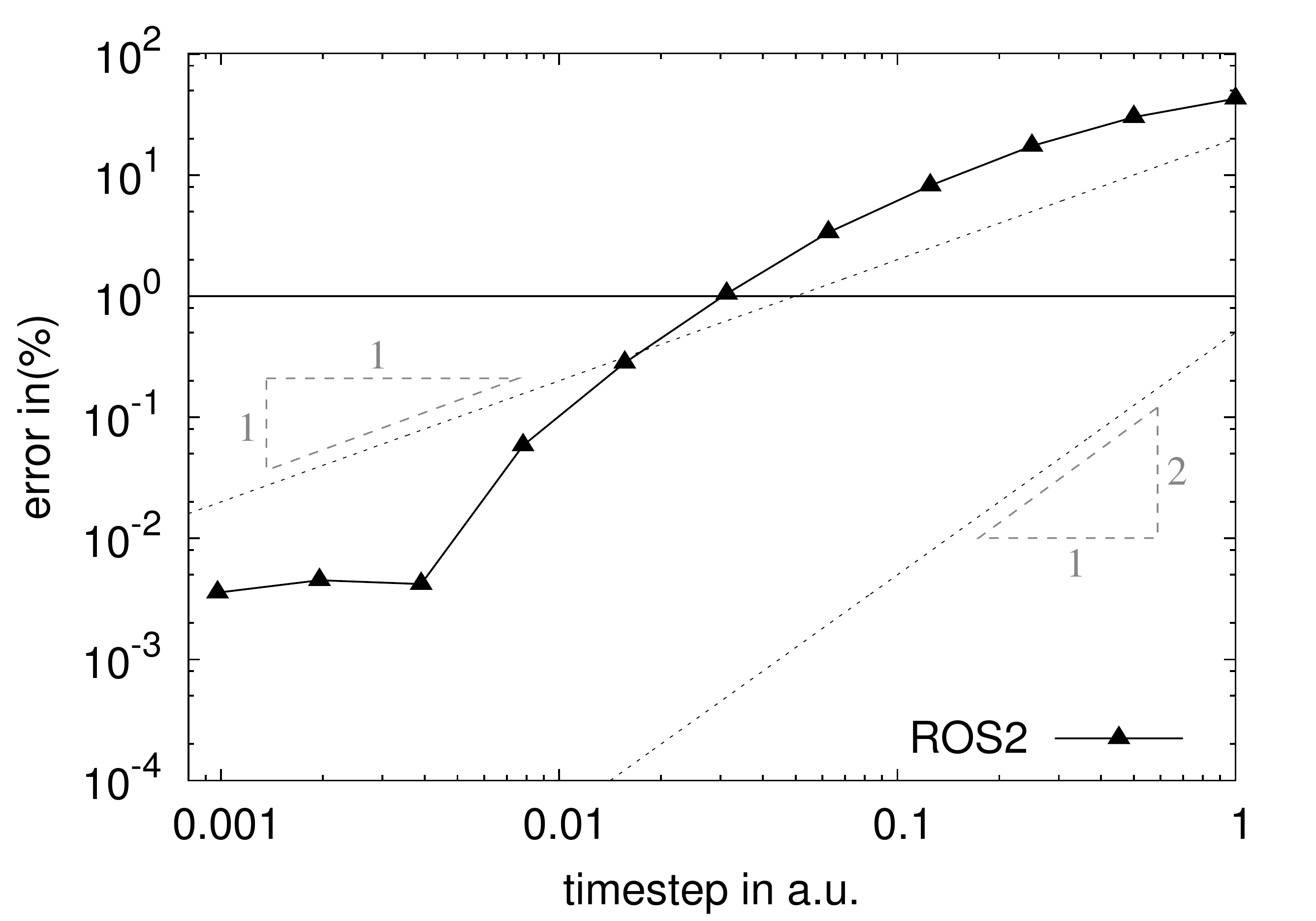} \\
\end{tabular}
\begin{center}
\begin{minipage}{0.9\textwidth}
\caption[short figure description]{
Rosenbrock ROS2 convexity splitting scheme. (left) Tip position over time for various timesteps together with the reference solution. The numbers indicate the used timesteps. (right) Deviation from the reference solution over timesteps.
\label{partfig3}
}
\end{minipage}
\end{center}
\end{figure}

\begin{figure}[htb]
\noindent \centering
\begin{tabular}{cc}
\includegraphics*[angle = -0, width = 0.45 \textwidth ]{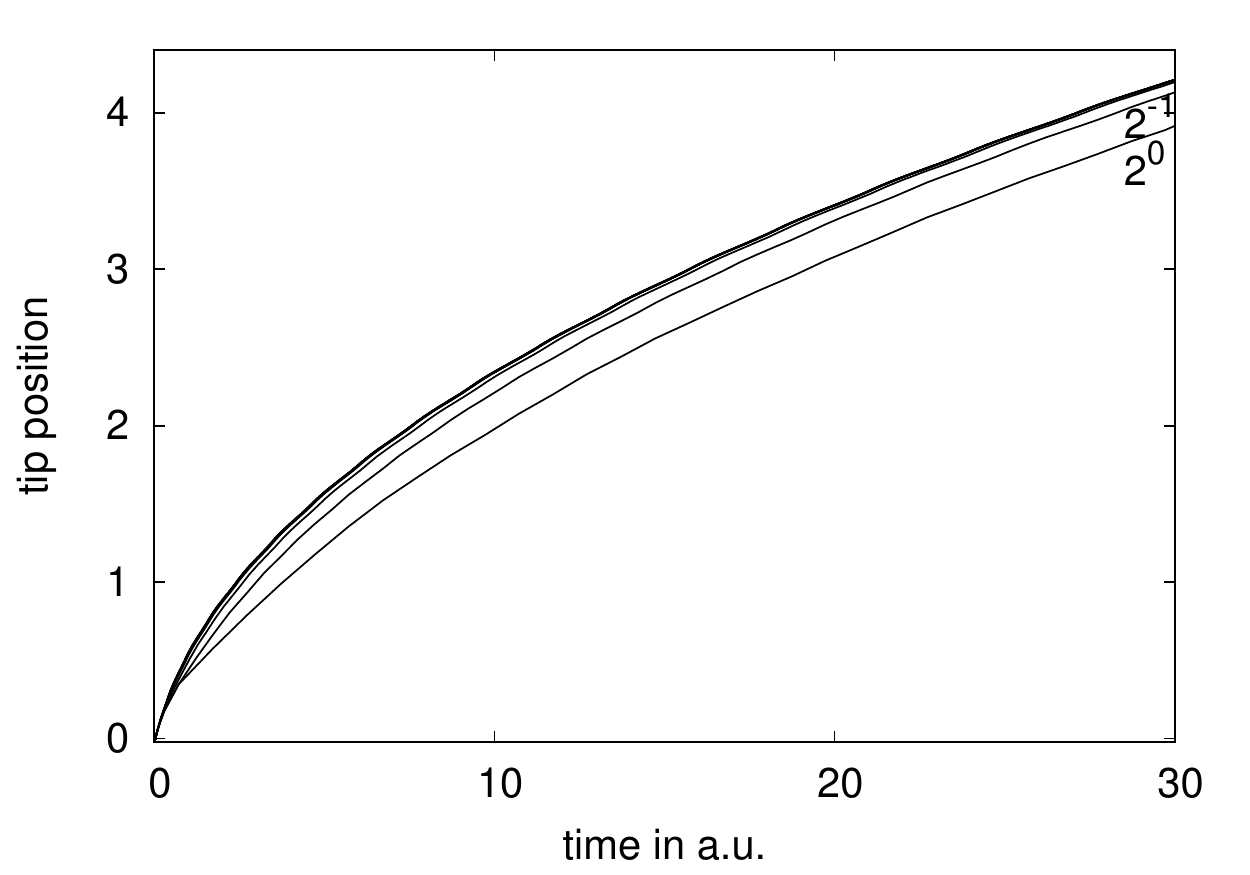} &
\includegraphics*[angle = -0, width = 0.45 \textwidth]{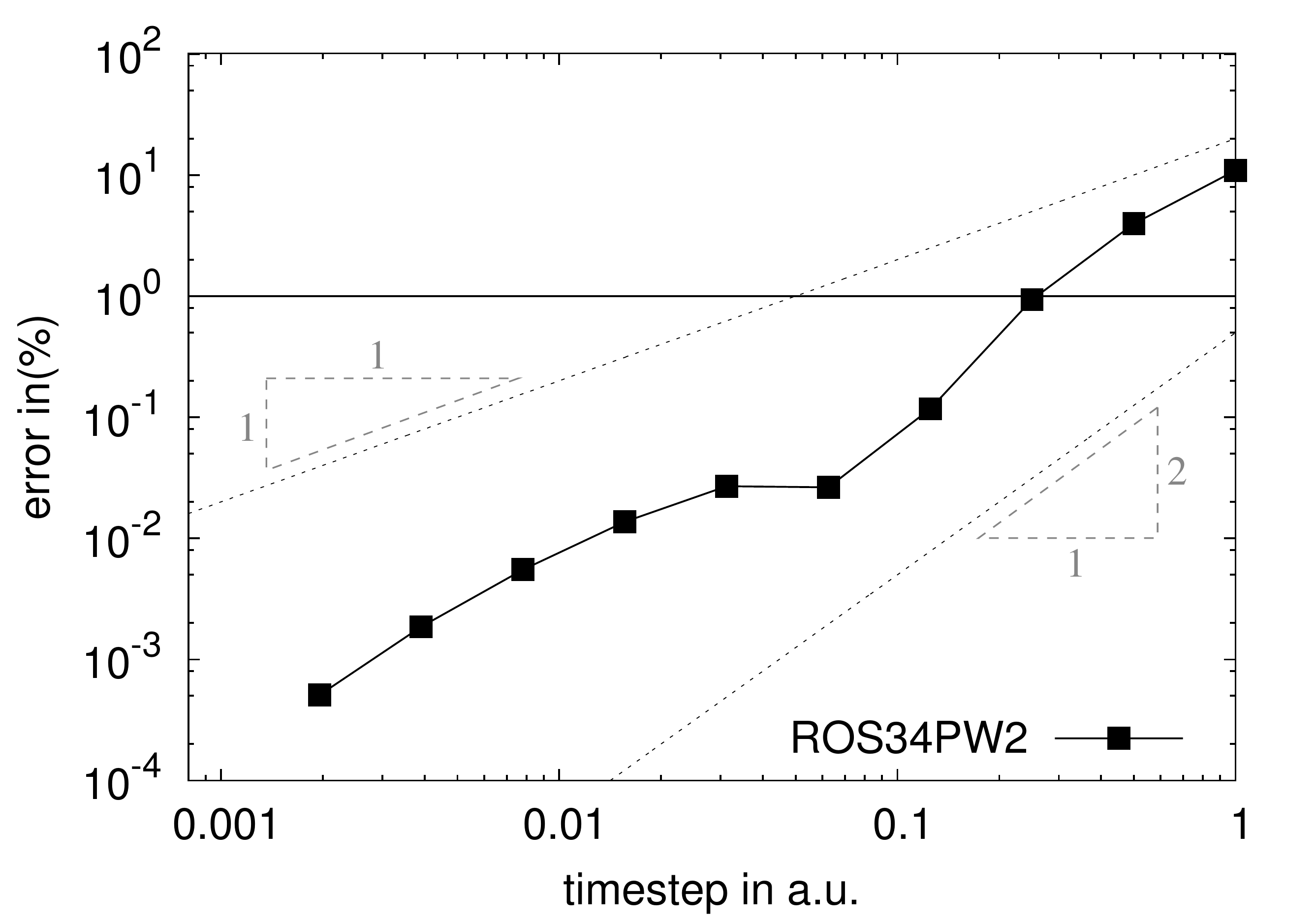} \\
\end{tabular}
\begin{center}
\begin{minipage}{0.9\textwidth}
\caption[short figure description]{
Rosenbrock ROS34WP2 convexity splitting scheme. (left) Tip position over time for various timesteps together with the reference solution. The numbers indicate the used timesteps. (right) Deviation from the reference solution over timesteps.
\label{partfig4}
}
\end{minipage}
\end{center}
\end{figure}

For the comparison with the experimental shapes in \cite{Naffoutietal_SA_2017} a numerical error within 1\% will be sufficient. This results from the uncertainties in the experimental measurements and the large modeling error. However, even if a numerical error within 1\% can be reached with the proposed scheme a much larger time step will be required to enable large scale simulations in three spatial dimensions. We demonstrate that this can be achieved using the Rosenbrock schemes. Figs. \ref{partfig3} and \ref{partfig4} show the results for ROS2 and ROS34WP2, respectively.

\begin{figure}[htb]
\noindent \centering
\begin{tabular}{c}
\includegraphics*[angle = -0, width = 0.45 \textwidth]{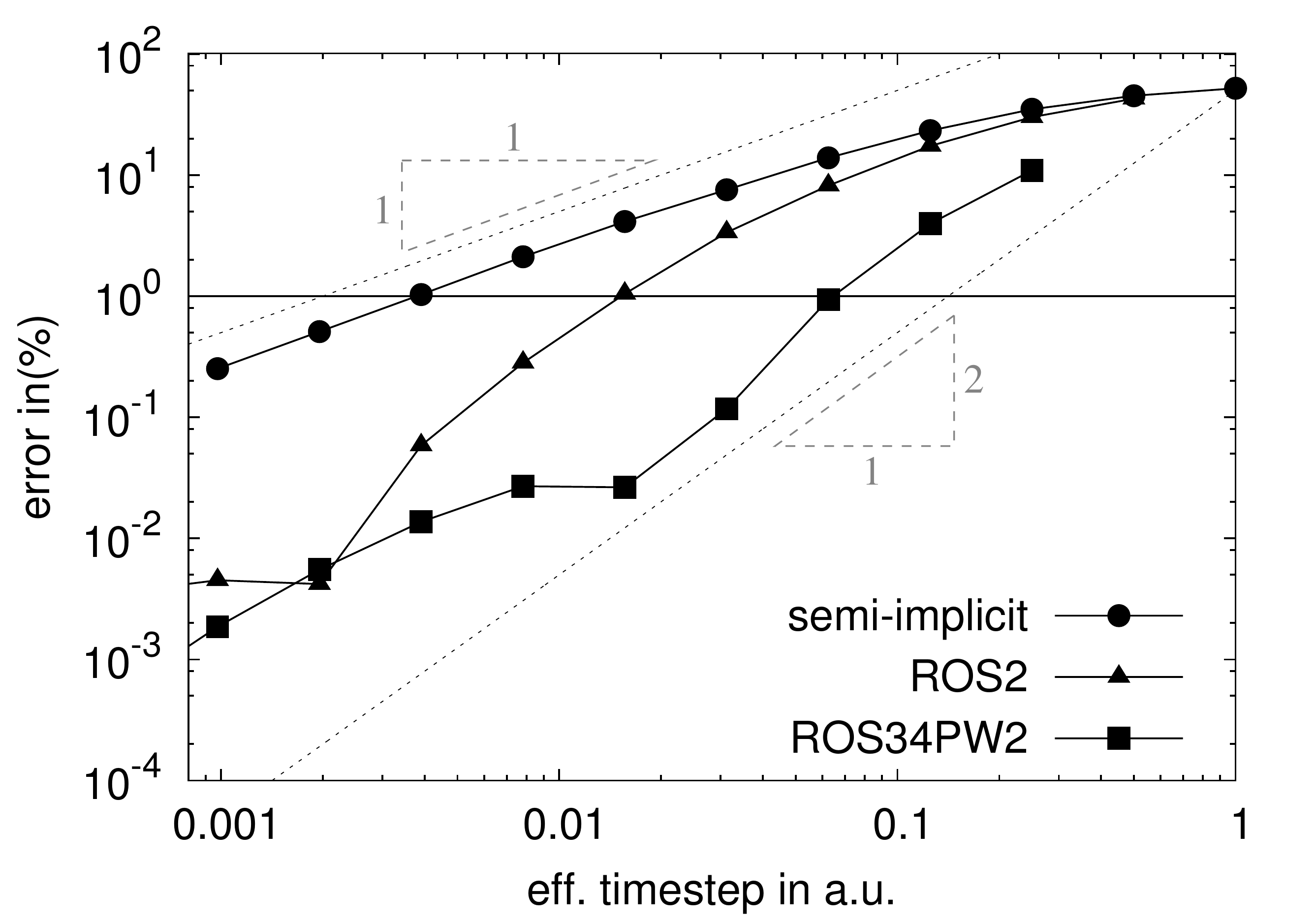} 
\end{tabular}
\begin{center}
\begin{minipage}{0.9\textwidth}
\caption[short figure description]{
Convergence w.r.t. effective numeric time step. Comparison of the semi-implicit convexity splitting and Rosenbrock convexity splitting schemes. 
\label{partfig5}
}
\end{minipage}
\end{center}
\end{figure}

The qualitative behavior of the reference solution can be recovered for $\tau < 10^{-1}$ and $\tau < 1$ for ROS2 and ROS34WP2, respectively. Quantitatively we obtain a solution with an error within 1\% for $\tau < 0.05 $ and $\tau < 0.25$. The results further indicate a better than first order convergence in $\tau$ in both cases. However, this improvement comes with an additional cost associated with the Rosenbrock schemes. The number of linear equations to be solved in each time step increases by a factor of two (ROS2), respectively four (ROS34WP2). We thus consider convergence with respect to an effective numerical time step $\tau_{\rm eff}=\tau / s$, with $s$ the number of steps in the considered Rosenbrock scheme. Fig. \ref{partfig5} shows the comparison of the three considered convexity splitting schemes. We observe first order convergence for the semi-implicit convexity splitting scheme and better than first order convergence for the Rosenbrock convexity splitting schemes. As the Jacobian is only approximated in our schemes we do not reach the theoretically predicted order of convergence of the Rosenbrock schemes. The ROS34WP2 scheme turns out to be the most efficient, it allows at least one order of magnitude larger timesteps than the semi-implicit convexity splitting scheme, without loss in accuracy. All schemes indicate the existence of a numerical upper bound for the timestep, a property which is also noticed for convexity splitting schemes where unconditional energy stability and unconditional solvability can be proven \cite{Chengetal_JCP_2008,Elseyetal_ESAIM_2013}.

We expect these results to be of general interest, especially in application where the long time behavior is concerned. In the next section we show large scale simulations, which would not be possible without the introduced convexity splitting Rosenbrock scheme. Other applications are found in phase field crystal simulations, where grain growth is considered and a Rosenbrock scheme already applied, see \cite{Backofenetal_AM_2014,Praetoriusetal_SIAMJSC_2015}. 

\section{Application}
\label{sec4}

\begin{figure}[htb]
\noindent \centering
\begin{tabular}{c}
\includegraphics*[angle = -0, width = 0.65 \textwidth]{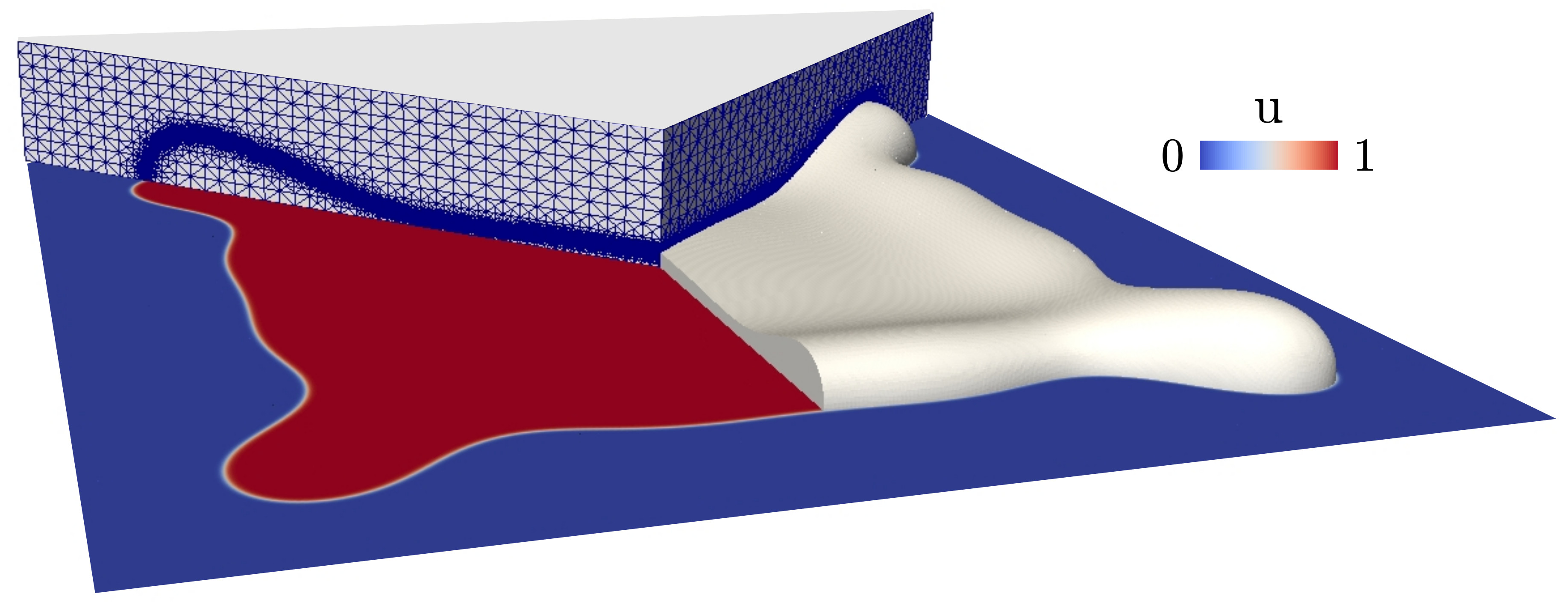} 
\end{tabular}
\begin{center}
\begin{minipage}{0.9\textwidth}
\caption[short figure description]{
Phase field variable, interface and adaptively refined mesh in order to ensure approximately 10 grid points across the interface.   
\label{partfig0}
}
\end{minipage}
\end{center}
\end{figure}

Fig. \ref{partfig0} shows the adaptively refined mesh and the phase field variable for a time step of the considered three dimensional simulations which are motivated by the experimentally observed nano-morphologies in \cite{Naffoutietal_SA_2017}.

\subsection{Numerical setting for nano-morphology simulation}
\begin{figure}[htb]
\noindent \centering
\begin{tabular}{cccccc}
A & B & C & D & E & F\\
\includegraphics*[angle = -0, width = 0.145 \textwidth]{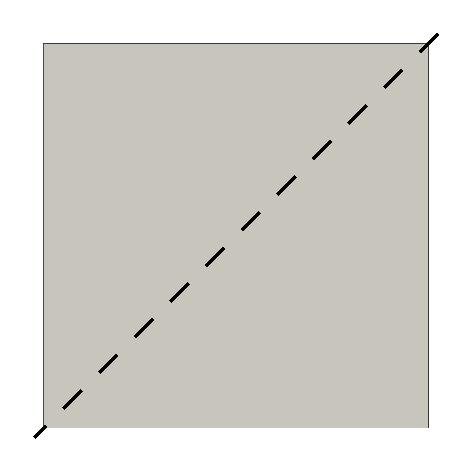} &
\includegraphics*[angle = -0, width = 0.145 \textwidth]{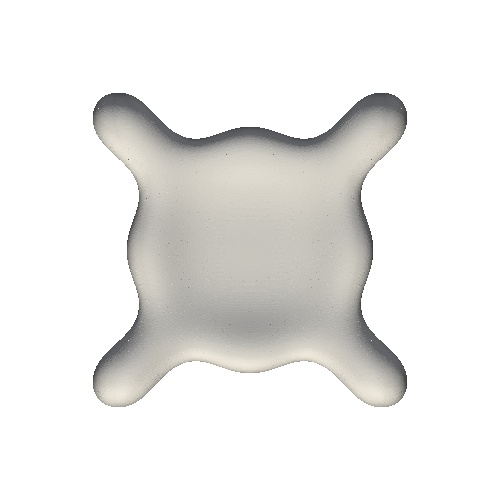} &
\includegraphics*[angle = -0, width = 0.145 \textwidth]{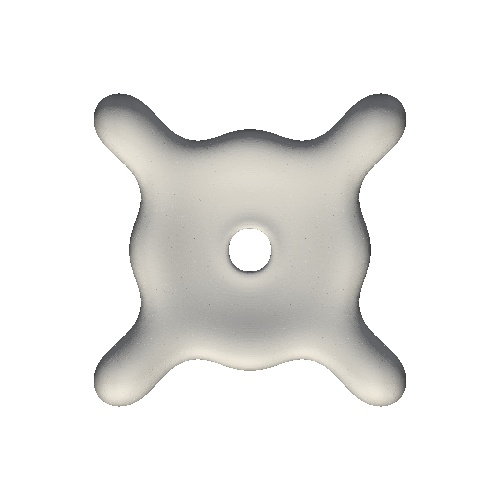} &
\includegraphics*[angle = -0, width = 0.145 \textwidth]{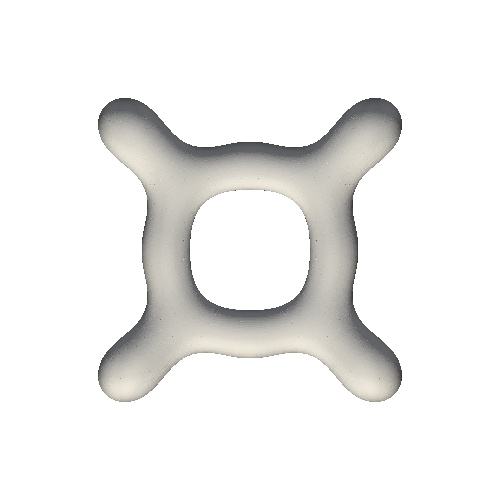} &
\includegraphics*[angle = -0, width = 0.145 \textwidth]{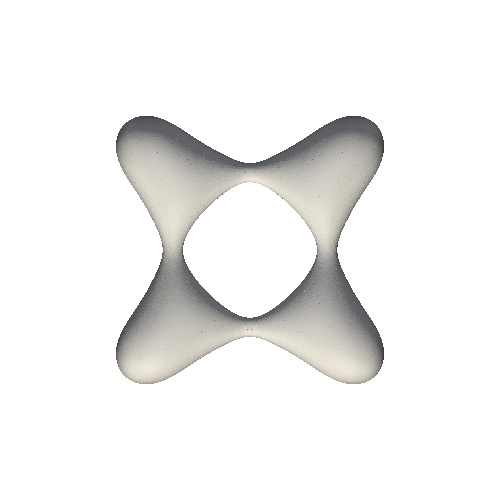} &
\includegraphics*[angle = -0, width = 0.145 \textwidth]{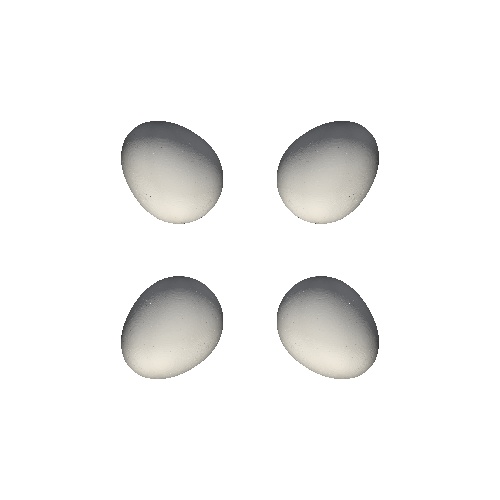} \\
\end{tabular}

\begin{tabular}{c}
\includegraphics*[angle = -0, width = 0.95 \textwidth]{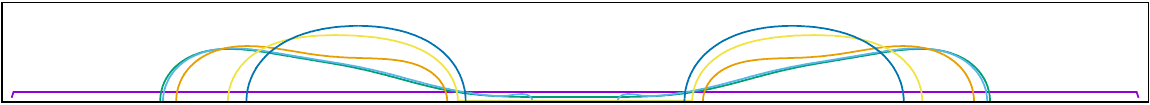} \\
\end{tabular}
\begin{center}
\begin{minipage}{0.9\textwidth}
\caption[short figure description]{
Dewetting of a square island. (top) Change in morphology from left to right and (bottom) height profile for different stages shown across the diagonal depicted in (A). The profiles corresponding to (A)-(F) evolve towards the center.
\label{partfig6}
}
\end{minipage}
\end{center}
\end{figure}

The ROS34WP2 convex splitting scheme is the method of choice for the large scale simulations in three spatial dimensions. Fig. \ref{partfig6} shows the results for the dewetting of a square island with aspect ratio height/length=$1/80$. (A) indicates the initial state. All sides retract but at the corners the retraction speed is smaller and fingers build up (B). The valley behind the tip eventually becomes so deep that the island breaks up and a hole is formed in the middle (C). This hole rapidly increases until it approaches the vicinity of the steps, thus, becoming square like (D). The bridges connecting the corner fingers become thinner (E) and break (F). The resulting four islands arrange as equidistant hemispherical dots (not shown).    

\begin{figure}[htb]
\noindent \centering
\begin{tabular}{c}
\includegraphics*[angle = -0, width = 0.65 \textwidth]{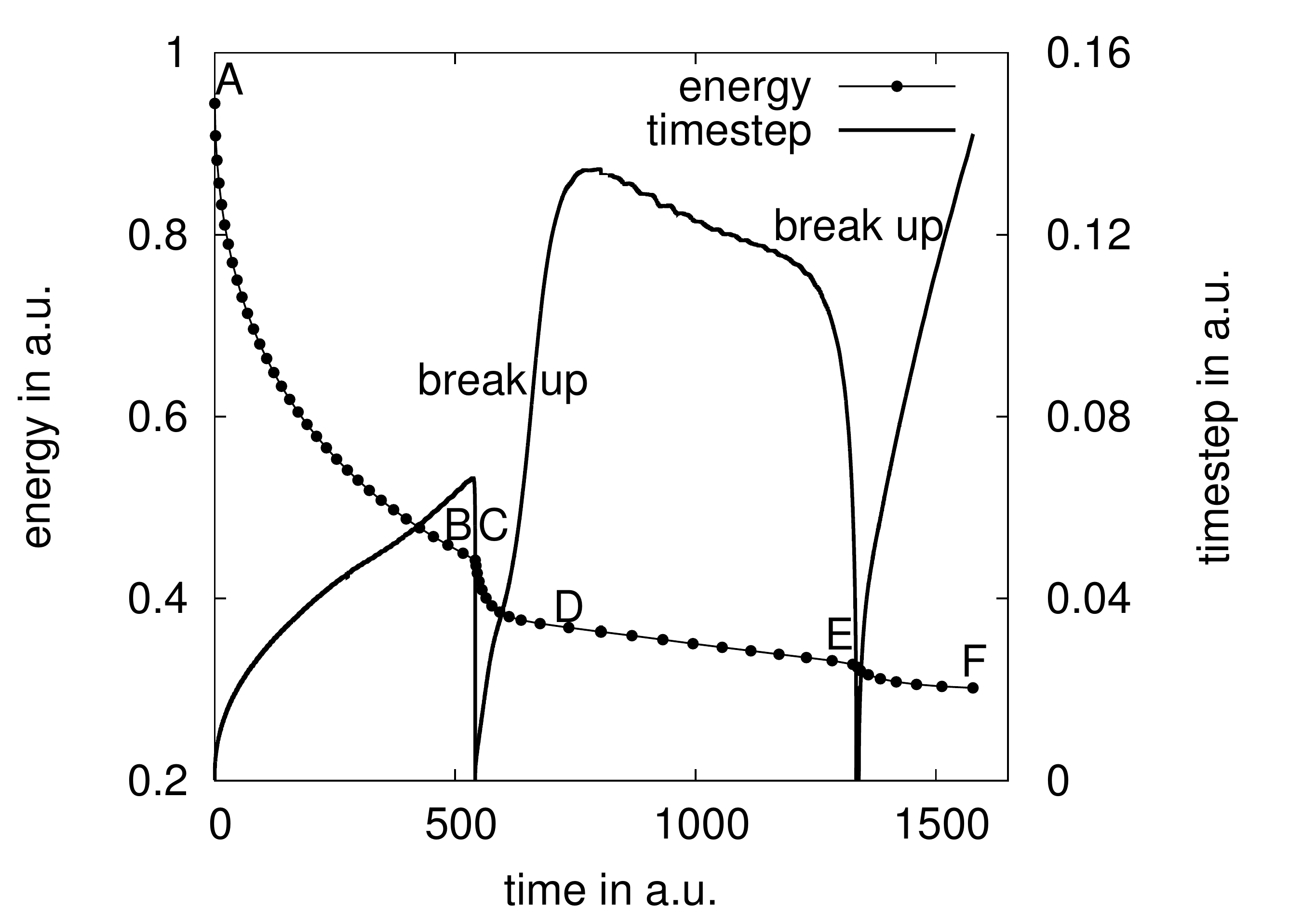} 
\end{tabular}
\begin{center}
\begin{minipage}{0.9\textwidth}
\caption[short figure description]{
  Dewetting of a square island. Energy decay and timestep evolution during the simulation. The snapshots (A)-(F) from Fig. \ref{partfig6} are labeled in the plot. The drastic reduction in energy associated with topological changes in the morphology can be seen, as well as the reduction in timestep associated with the increased dynamics of such events. For smooth morphology evolutions large timesteps, close to the maximal numerical timestep to ensure the required accuracy are chosen.
\label{partfig7}
}
\end{minipage}
\end{center}
\end{figure}

To enable these simulations we exploit the symmetry of the system and only calculate a quarter of the domain. We further make use of an additional advantage of the Rosenbrock scheme. It allows to compute a lower order approximation of the solution without much additional cost \cite{Lan99},
\begin{eqnarray*}
\uv^{n+1}_{\rm low}=\uv^{n}+\sum_{i = 1}^{s} \mv_{i} \uv^{n}_i \,
\end{eqnarray*}
This allows for a proper definition of time errors $e^{n+1}= ||\uv^{n+1}-\uv^{n+1}_{\rm low} ||$ which can be used to adapt the timestep \cite{Johnetal_CMAME_2010}. The next timestep is e.g. controlled by a PI-controller \cite{Lan99}, 
\begin{eqnarray*}
\tau^{n+1} = \rho \frac{(\tau^{n})^2}{\tau^{n-1}} \left( \frac{e_{\rm tol} e^{n}}{(e^{n+1})^2}\right)^{1/p},
\end{eqnarray*} 
where $e_{\rm tol}$ is a prescribed error bound, $\rho \in (0,1]$ a relaxation factor and $p$ the order of the Rosenbrock method. In the following we use $e_{\rm tol} = 4 \cdot 10^{-3}$, $\rho = 0.95$ and $p = 3$.

Fig. \ref{partfig7} shows the reduction of the energy and the considered timestep for the simulation in Fig. \ref{partfig6}. The various stages of the evolution indicated by (A) - (F) are shown and demonstrate the relation between a drastic reduction in timestep and topological changes in the morphology and the large timesteps used for smooth morphology changes. The simulations are run on a parallel environment with 384 cores, using domain decomposition. 

\subsection{Model extension}

Even if quantitative comparisons between the experimentally observed and the computed nano-morphologies are already possible, see Fig. 4 in \cite{Naffoutietal_SA_2017}, improvements in the considered model are needed
to further reduce the discrepancies. Besides the process condition this includes the incorporation of vapor-substrate and film-substrate interfacial energies and anisotropy. 

\subsubsection{Wetting angle}
\begin{figure}[htb]
\noindent \centering
\begin{tabular}{c}
\includegraphics*[angle = -0, width = 0.5 \textwidth]{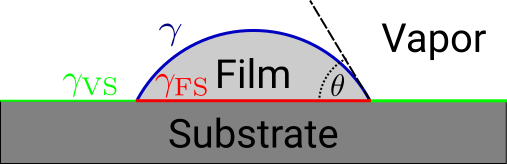} 
\end{tabular}
\begin{center}
\begin{minipage}{0.9\textwidth}
\caption[short figure description]{
Young's law. The relation between (isotropic) film-vapor interface energy ($\gamma$), film-substrate interface energy ($\gamma_{\rm FS}$) and the vapor-substrate interface energy ($\gamma_{\rm VS}$) ldetermines to a unique wetting angle $\theta$. 
\label{fig:YoungsLaw}
}
\end{minipage}
\end{center}
\end{figure}

Following typical approaches for contact problems in fluid dynamics  \cite{Jacqmin_JFM_2000,Villanuevaetal_IJMPF_2006,Alandetal_CMES_2010} we introduce the substrate energy
\begin{eqnarray*}
{\mathcal{E}}_{sub}[u] = \int_\Omega \frac{1}{\xi} B(v) \left(\frac{1}{2}(\gamma_{\rm VS} + \gamma_{\rm FS}) - \frac{-4 u^3 + 6 u^2 - 1}{2}(\gamma_{\rm VS} - \gamma_{\rm FS}) \right) \; d \mathbf{x}
\end{eqnarray*}
with vapor-substrate and film-substrate energy densities, $\gamma_{\rm VS}$ and $\gamma_{\rm FS}$, respectively. We define $v(z) = \frac{1}{2}(1 - \tanh(\frac{3 z}{\xi}))$, with $z$ the height above the substrate and $\xi > 0$ a small parameter. $\frac{1}{\xi} B(v)$ is thus an approximation for a delta function used to consider the boundary condition at the substrate. The cubic polynomial in $u$ ensures the substrate energy density to be equal to $\gamma_{\rm VS}$ for $u = 0$ and to be equal to $\gamma_{\rm FS}$ for $u = 1$, as well as it's derivative to be zero for $u = 0$ and $u = 1$. This energy now has to be added to $\mathcal{E}$ in the derivation of the evolution equations, which leads instead of eq. \eqref{eq1} to
\begin{eqnarray}
\label{eq1sub}
\partial_t u &=& \nabla \cdot \mathbf{j}, \qquad
\mathbf{j} = \frac{1}{\epsilon} M(u) \nabla \mu, \\
\label{eq3sub} 
g(u)\mu &=& \frac{1}{\epsilon} B^\prime(u) - \epsilon \Delta u + \frac{1}{\xi} B(v) 6 u (u -1) (\gamma_{\rm VS} - \gamma_{\rm FS}),
\end{eqnarray}
in $\Omega \times(0,\infty)$. The initial and boundary conditions remain. Following \cite{Lietal_CMS_2009}, the asymptotic limit $\xi \to 0$ leads to eq. \eqref{eq1} with boundary condition $\mathbf{n} \cdot \nabla \mu = 0$ and $\epsilon \mathbf{n} \cdot \nabla u = 6 u ( u - 1) (\gamma_{\rm VS} - \gamma_{\rm FS})$. Using Young's law,  $\gamma_{\rm VS} - \gamma_{\rm FS} = \gamma \cos (\theta)$, with equilibrium contact angle $\theta$ and film-vapor energy density $\gamma$, which in our case is constant and equal to one, this is consistent with the treatment in \cite{Jiangetal_AM_2012}. A similar approach has recently been proposed in \cite{Dziwniketal_Nonl_2017}. With $\theta = 90^\circ$ or equivalently $\gamma_{\rm VS} = \gamma_{\rm FS}$ we obtain our previous model. The sharp interface limit $\epsilon \to 0$ for the above treatment of the triple junction is considered in \cite{NovickCohen_PhysD_2000,Yueetal_JFM_2010} and leads to the classical Young's law, see Fig.~\ref{fig:YoungsLaw}. 

The formulation in eqs. \eqref{eq1sub} - \eqref{eq3sub} allows to use the proposed convexity splitting approach. We only need to modify eq. \eqref{eq4}, which now read
\begin{eqnarray}
\label{eq4sub}
B_c(u)&=&B(u) + B(v) (6 u(u -1) (\gamma_{\rm VS} - \gamma_{\rm FS}) + \alpha  \left(u-\frac{1}{2}\right)^2, \\
\label{eq5sub} 
B_e(u)&=&\alpha \left(u-\frac{1}{2}\right)^2,
\end{eqnarray}
where we have set $\xi = \epsilon$. To ensure convexity $\alpha$ now depends on $\left(\gamma_{\rm VS}-\gamma_{\rm FS}\right)$ and reads
\begin{eqnarray}
\alpha \geq \alpha_0 +\frac{B(v)^2}{12 \epsilon} (\gamma_{\rm VS} - \gamma_{\rm FS}),
\end{eqnarray}
with $\alpha_0 \geq 9$ as above.  The proposed formulation in eqs. \eqref{eq1sub} - \eqref{eq3sub} furthermore has the advantage to circumvent the definition of a contact angle, which becomes less meaningful if anisotropies are considered. We first analyze the effect in the isotropic case.

\begin{figure}[htb]
\noindent \centering
\begin{tabular}{cccccc}
A & B & D & D & E & F\\
\includegraphics*[angle = -0, width = 0.145 \textwidth]{./40_0_0_nice} &
\includegraphics*[angle = -0, width = 0.145 \textwidth]{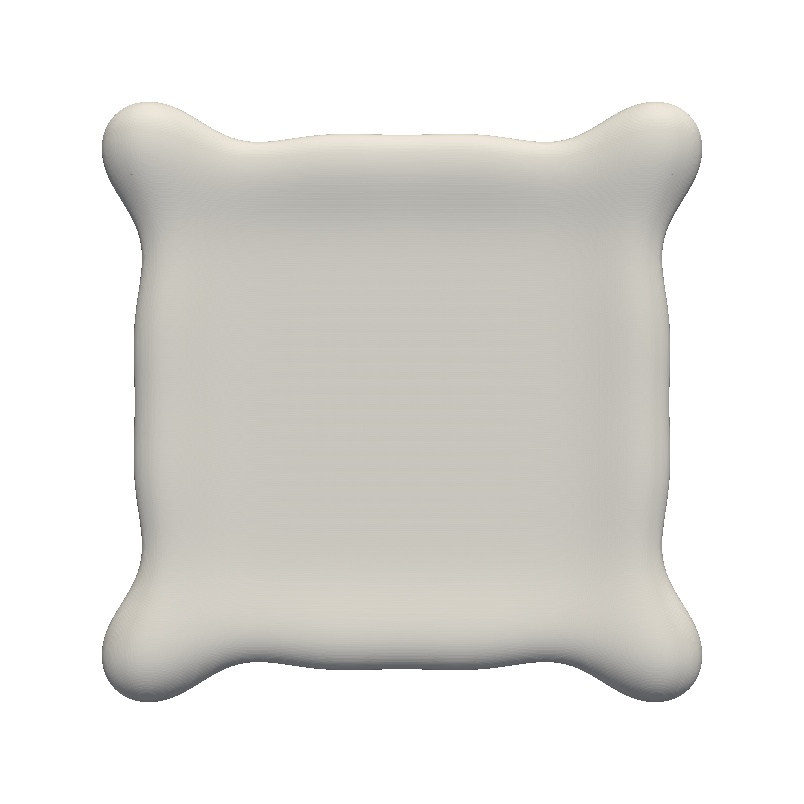} &
\includegraphics*[angle = -0, width = 0.145 \textwidth]{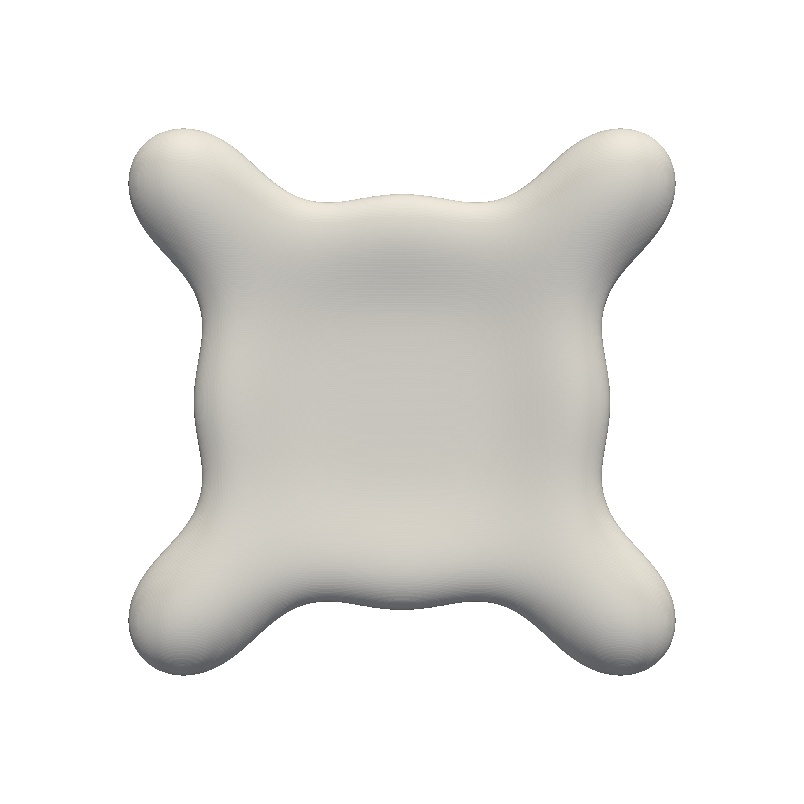} &
\includegraphics*[angle = -0, width = 0.145 \textwidth]{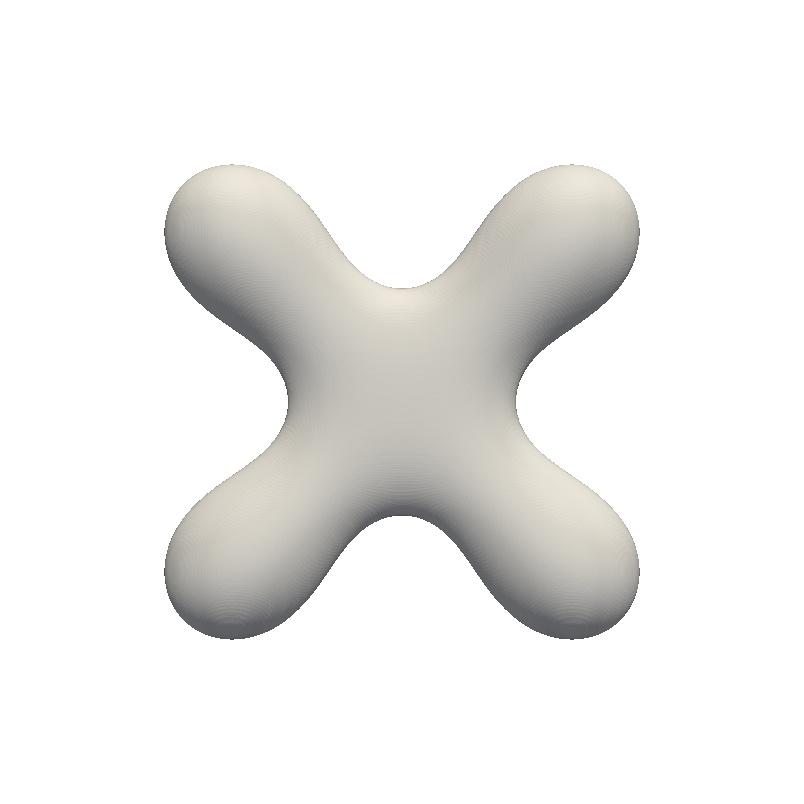} &
\includegraphics*[angle = -0, width = 0.145 \textwidth]{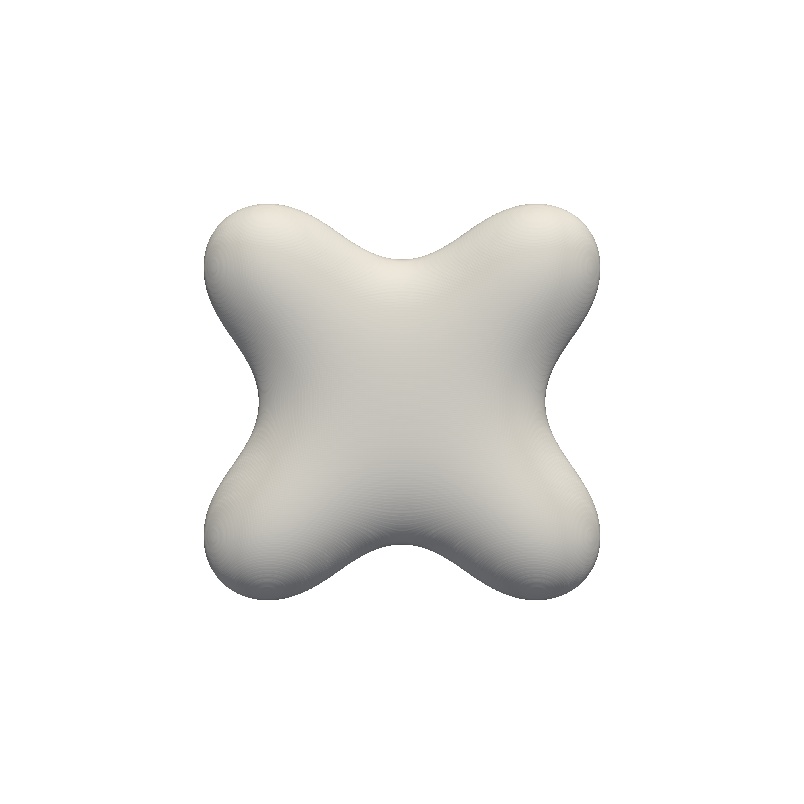} &
\includegraphics*[angle = -0, width = 0.145 \textwidth]{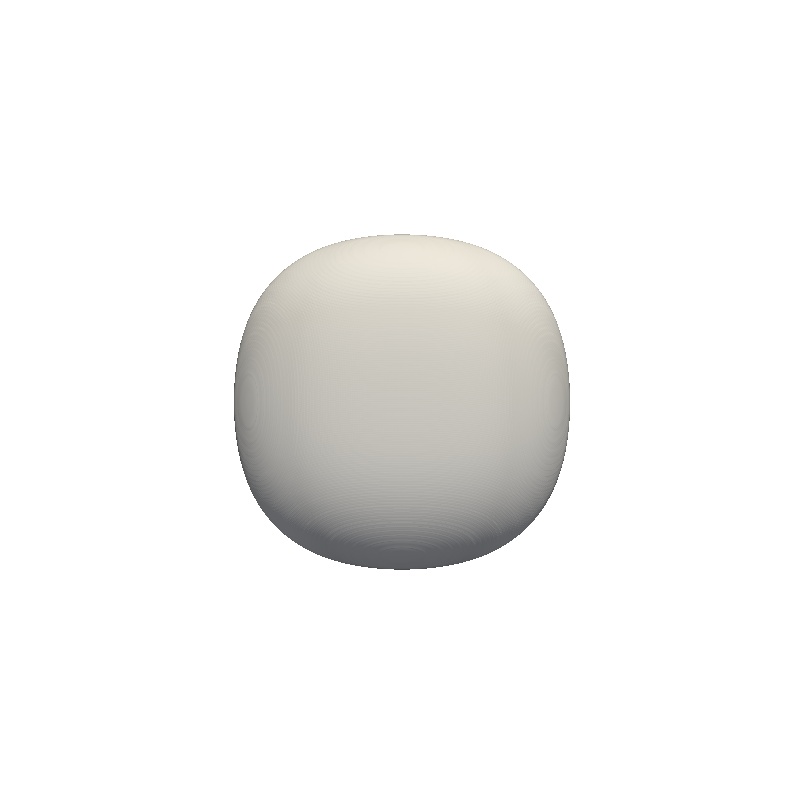} \\
\end{tabular}
\begin{tabular}{c}
\includegraphics*[angle = -0, width = 0.95 \textwidth]{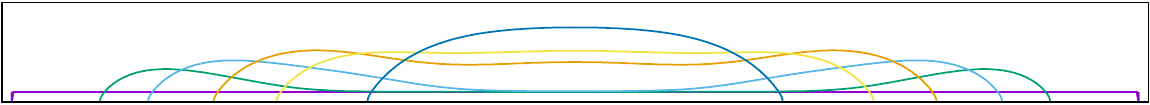} \\
\end{tabular}
\begin{center}
\begin{minipage}{0.9\textwidth}
\caption[short figure description]{
Dewetting of a square island. (top) Change in morphology from left to right and (bottom) height profile for different stages shown across the diagonal depicted in (A). The profiles corresponding to (A)-(F) evolve towards the center. The corresponding times of the snap-shots are (A,...,F) = (0,219, 670, 1610, 3422, 4802).
\label{partfig10}
}
\end{minipage}
\end{center}
\end{figure}

\begin{figure}[htb]
\noindent \centering
\begin{tabular}{cccccc}
A & B & C & D & E & F\\
\includegraphics*[angle = -0, width = 0.145 \textwidth]{./40_0_0_nice} &
\includegraphics*[angle = -0, width = 0.145 \textwidth]{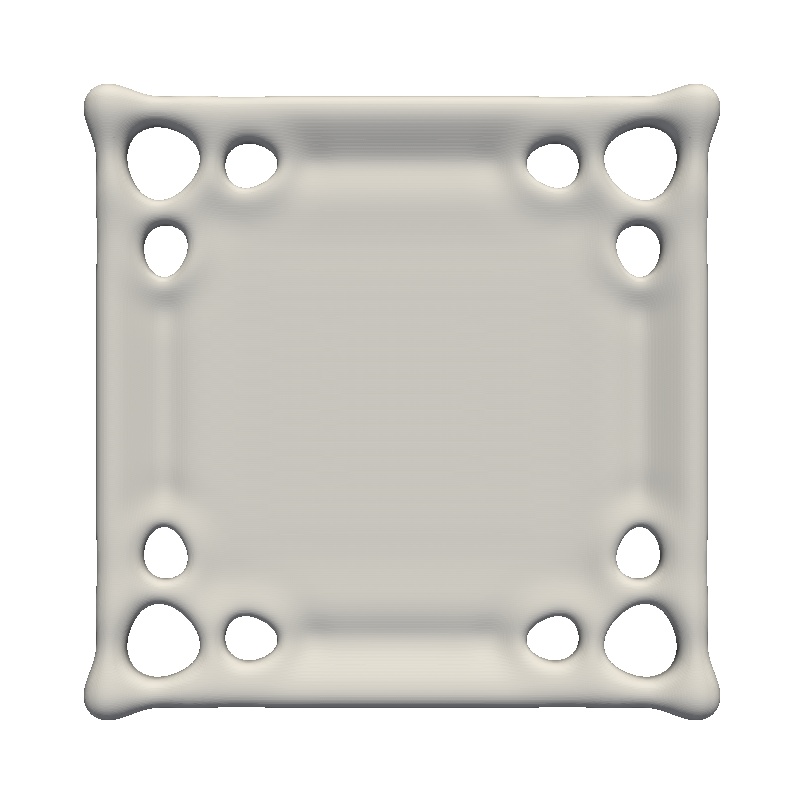} &
\includegraphics*[angle = -0, width = 0.145 \textwidth]{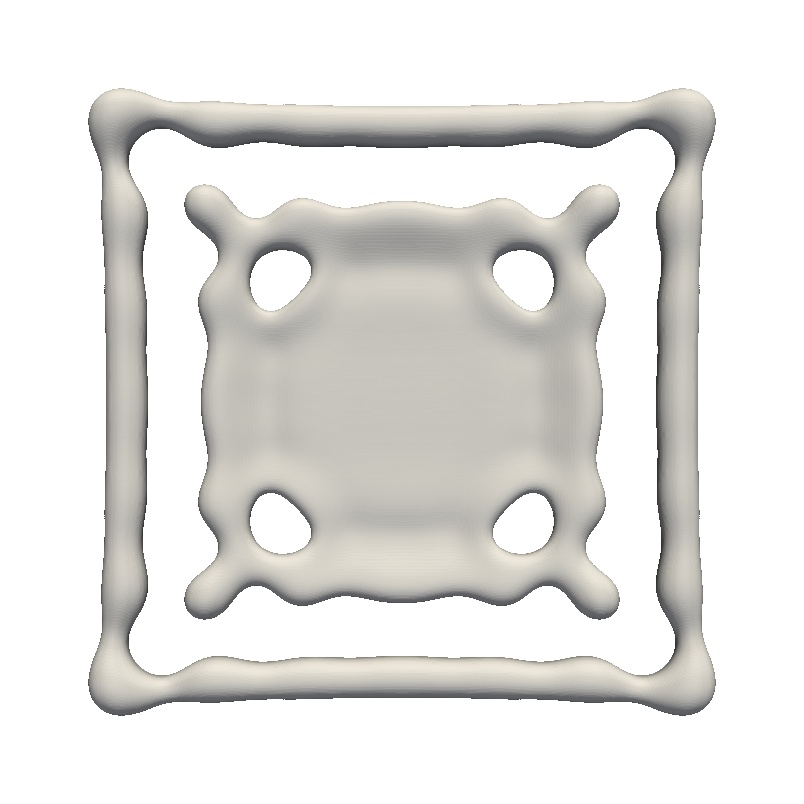} &
\includegraphics*[angle = -0, width = 0.145 \textwidth]{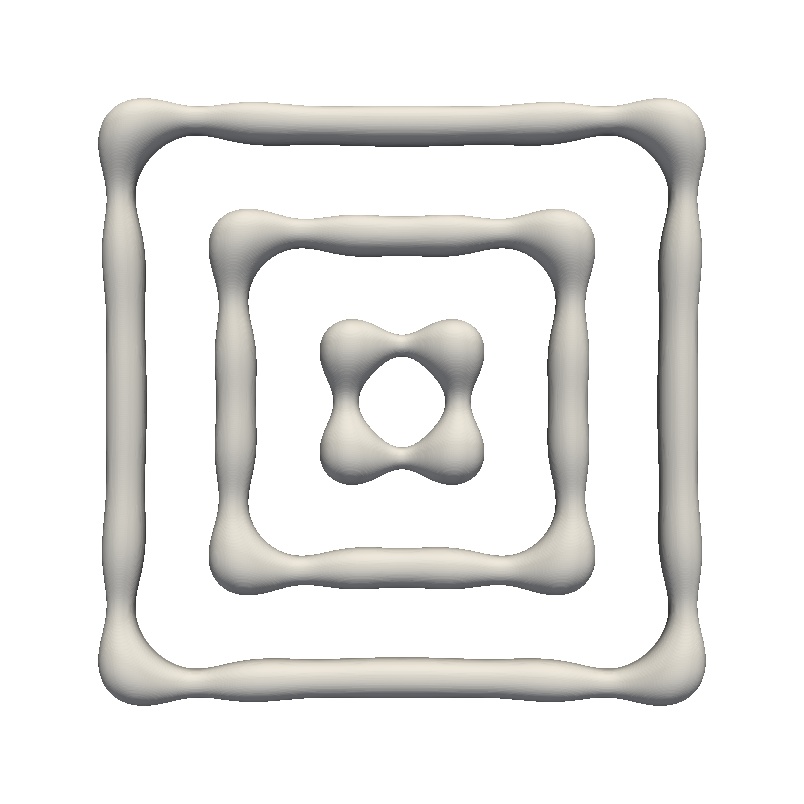} &
\includegraphics*[angle = -0, width = 0.145 \textwidth]{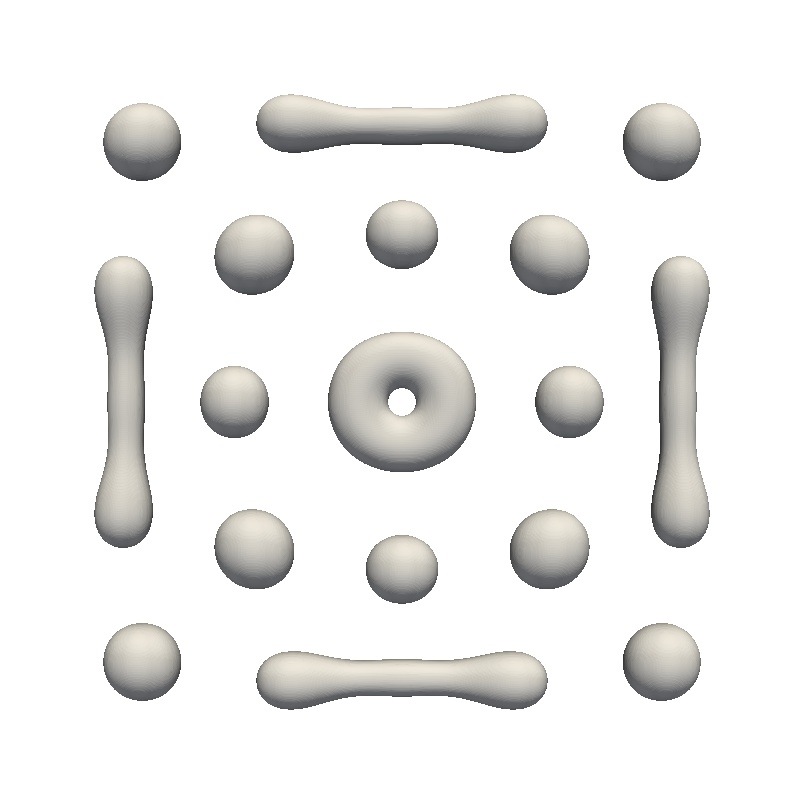} &
\includegraphics*[angle = -0, width = 0.145 \textwidth]{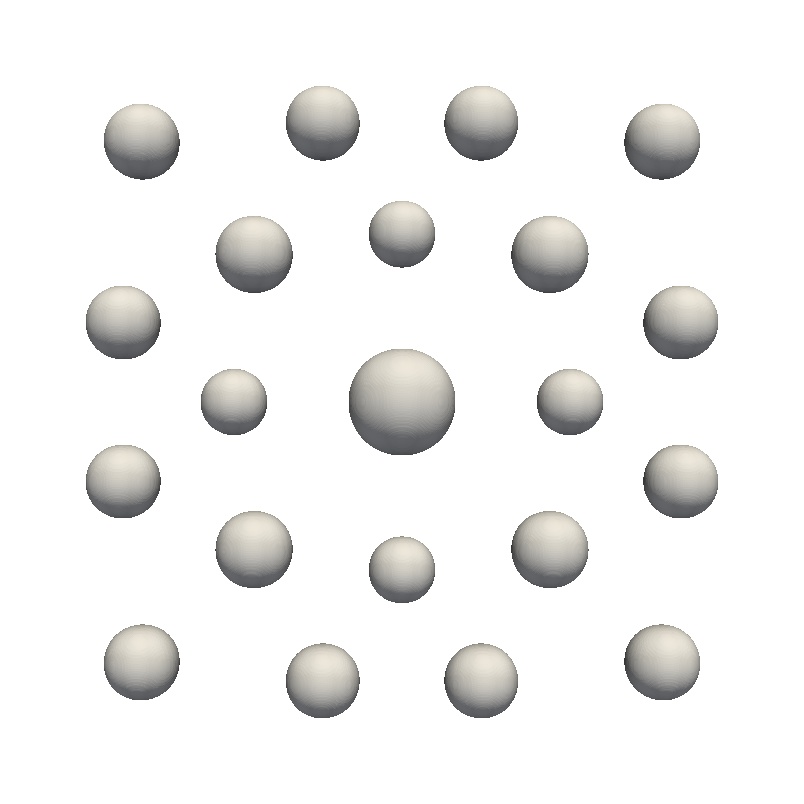} \\
\end{tabular}

\begin{tabular}{c}
\includegraphics*[angle = -0, width = 0.95 \textwidth]{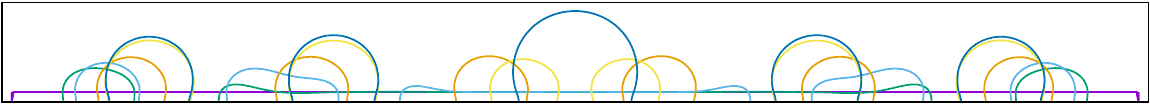} \\
\end{tabular}
\begin{center}
\begin{minipage}{0.9\textwidth}
\caption[short figure description]{
Dewetting of a square island. (top) Change in morphology from left to right and (bottom) height profile for different stages shown across the diagonal depicted in (A). The profiles corresponding to (A)-(F) evolve towards the center. The corresponding times of the snap-shots are (A,...,F) = (0,22, 40, 123, 274, 1123).
\label{partfig11}
}
\end{minipage}
\end{center}
\end{figure}

We consider two scenarios $\nicefrac{(\gamma_{\rm VS} - \gamma_{\rm FS})}{\gamma}=0.5$ and $\nicefrac{(\gamma_{\rm VS} - \gamma_{\rm FS})}{\gamma}=-0.5$, corresponding to wetting angles $\theta= 60^\circ$ and $\theta= 120^\circ$, respectively. Fig.~\ref{partfig10} shows the evolution for the first case. All sides of the initial square (A) retract. The retraction speed at the corners is smaller which leads to the formation of fingers at the corners (B). Due to the smaller wetting angle, the hill and valley behind the retracting front is elongated and not as pronounced (C). Thus, the breaking of the film due to hole formation is suppressed. Instead, the evolving fingers at the corners lead to a cross-shape (D). This shape then becomes more and more compact and evolves to a singular drop (E-F).  The final shape is a sphere, which is cut by the substrate to fullfil the volume constraint and the equilibrium wetting angle (not shown). The second case is shown in Fig.~\ref{partfig11}, where a completely different effect is observed.

The hill behind the retracting step grows faster and the valley becomes deeper. Thus, the film touches the substrate earlier and forms holes (B). At the corners the valley is thiner due to the effect of both adjacent steps and thus the hole formation starts earlier. A chain of holes  starts to form in the primary valley parallel to the sides of the initial structure. This holes join and separate an outer set of material lines from a inner square patch (C). The outer lines form a square and decompose similar to a Rayleigh instability and form a line of dots (D-F). The inner square patch resembles the initial square on a smaller scale. The remaining patch in the middle, is too small to further decompose and collapse to a single drop (D-F).  At the end a regular array of dots is achieved. Every dot is a sphere cut by the substrate with a wetting angle, $\Theta = 120^\circ$.   

We should mentioned that the observed topological changes are quite sensitive to the interface width of the phase field and the width of the surface delta function. Smaller $\epsilon$ lead to a delay of the touching of the interface with the substrate. Thus, the hole formation occurs later. However, the general trend is independent on the modeling details. A reduced contact angle $\theta < 90^\circ$ leads to more compact shapes, while an increased contact angle $\theta > 90^\circ$ enhances further splitting. 

In our modeling approach the interaction between film and substrate is smeared out by construction. Due to the boundary condition the isolines of the phase field are still forced to touch the substrate with $90^\circ$ on a smaller length scale than the interface width.  Thus the wetting angle has to be defined by extrapolating the shape of the film towards the substrate. In the equilibrium state, a sphere may be fitted to the island shape. The angle of the sphere at the substrate then defines a proper wetting angle.   

\subsubsection{Anisotropy}

We refer to
\cite{Torabietal_PRSA_2009} for a treatment of weak and strong anisotropies in the film-vapor interfacial energy density in the context of phase field approximations for surface diffusion. The transition between weak and strong anisotropies occurs at a convex-to-concave transition in the $1/ \gamma$ plot \cite{Sekerka_JCG_2005}, where $\gamma = \gamma(\boldsymbol{\nu})$, with $\boldsymbol{\nu} = \nabla u / |\nabla u|$ the normal to the film-vapor interface. For weak anisotropies the equations without the consideration of substrate interfacial energies now read
\begin{eqnarray}
\label{eq1a}
\partial_t u &=& \nabla \cdot \mathbf{j}, \qquad \mathbf{j} = \frac{1}{\epsilon} M(u) \nabla \mu, \\
\label{eq3a} 
g(u)\mu &=& \frac{1}{\epsilon} \gamma(\boldsymbol{\nu}) B^\prime(u) - \epsilon \nabla \cdot (\gamma(\boldsymbol{\nu}) \nabla u) - \epsilon \nabla \cdot (|\nabla u| ^2 \nabla_{\nabla u} \gamma(\boldsymbol{\nu)}).
\end{eqnarray}
For strong anisotropies these equations become ill-posed and require a regularization, see \cite{Gurtinetal_ARMA_2002,Friedetal_APM_2004}. Following \cite{Torabietal_PRSA_2009} a Willmore regularization is added to penalize large curvatures with penalization parameter $\beta$. The equations now read
\begin{eqnarray}
\label{eq1b}
\partial_t u &=& \nabla \cdot \mathbf{j}, \qquad
\mathbf{j} = \frac{1}{\epsilon} M(u) \nabla \mu, \\
\label{eq3b} 
g(u)\mu &=& \frac{1}{\epsilon} \gamma(\boldsymbol{\nu}) B^\prime(u) - \epsilon \nabla \cdot (\gamma(\boldsymbol{\nu}) \nabla u) - \epsilon \nabla \cdot (|\nabla u| ^2 \nabla_{\nabla u} \gamma(\boldsymbol{\nu)}) \nonumber \\
&& + \beta (\frac{1}{\epsilon^2} B^{\prime\prime}(u) \kappa - \Delta \kappa), \\
\kappa &=&  \frac{1}{\epsilon} B^\prime(u) - \epsilon \Delta u.
\end{eqnarray}
Both models follow from the energy
\begin{eqnarray*}
{\mathcal{E}}(u) = \int_\Omega \gamma(\boldsymbol{\nu}) \left( \frac{\epsilon}{2} |\nabla u|^2 + \frac{1}{\epsilon} B(u) \right) \; d \mathbf{x} + \frac{\beta}{2 \epsilon} \int_\Omega (- \epsilon \Delta u + \frac{1}{\epsilon} B^\prime(u))^2 \; d \mathbf{x}
\end{eqnarray*}
with $\beta = 0$ and $\beta > 0$, respectively.. However, in both cases the asymptotic result that near interfaces $\frac{\epsilon}{2} |\nabla u| ^2 \sim \frac{1}{\epsilon} B(u)$ is used to simplify the expression, see \cite{Lietal_CCP_2009}. The previously considered eq. \eqref{eq1} are obtained for $\gamma(\boldsymbol{\nu}) = 1$. A functional form for the surface energy density $\gamma(\boldsymbol{\nu})$, which is suitable for a specific material can be obtained from the approach proposed in \cite{Salvalaglioetal_CGD_2014}. Material specific simulations with strong anisotropies can be found in \cite{Salvalaglioetal_ACSAMI_2015,Bergamaschinietal_APX_2016,Salvalaglioetal_ASS_2017,Salvalaglioetal_NRL_2017}.   

Large scale simulations for these models which allow a detailed investigation of solid-state dewetting are still work in progress. The same is true for a combination of substrate interfacial energies and anisotropies. All current results in this direction only consider two-dimensional models, see \cite{Dziwniketal_Nonl_2017,Baoetal_JCP_2017}.

\section{Conclusions}
\label{sec5}

The advantages of convexity splitting schemes, which might be unconditionally energy stable, unconditionally solvable and optimally convergent in the energy norm, come with a reduction in accuracy. For simulations within a given error bound a maximal numerical timestep exist. This maximal numerical timestep might not be much larger than in classical schemes without convexity splitting. Even if this disadvantage has been pointed out by several people \cite{Chengetal_JCP_2008,Christliebetal_CMS_2013,Elseyetal_ESAIM_2013,Glasneretal_JCP_2016} many examples exist where this fact is not respected and convexity splitting schemes are used with unrealistic large timesteps. We here propose a convexity splitting scheme with increased accuracy for a phase field model for surface diffusion \cite{Raetzetal_JCP_2006}. The convexity splitting idea is combined with a Rosenbrock time stepping scheme. Through various approximations we make large scale simulations in three spatial dimensions tractable.
We numerically demonstrate the accuracy of the method on an example in two spatial dimensions. The considered convexity splitting Rosenbrock scheme ROS34RW2 allows for at least one order of magnitude larger time steps than the semi-implicit convexity splitting scheme without loss in accuracy. We demonstrate the possibilities the proposed scheme offers for exploring  the physical phenomena behind solid-state dewetting in three spatial dimensions. 

\section*{Acknowledgments}
AV and RB acknowledge support by the German Research Foundation (DFG) under grant SPP1959, MS acknowledges support by the Postdoctoral Research Fellowship awarded by the Alexander von Humboldt Foundation and SW acknowledges support by the DRESDEN Fellowship program. We further acknowledge computing resources provided at J\"ulich Supercomputing Center under grant HDR06.

\bibliography{lit}
\bibliographystyle{plain}
\end{document}